\documentclass[twoside]{article}
\newcommand{\heuteIst}{July 11, 2002}

\usepackage{amsmath,amsthm,latexsym,amscd,amssymb}
\newcommand{\href[2]}{#2}

\usepackage{amsthm}

\swapnumbers
\theoremstyle{plain}
\newtheorem{theorem}{Theorem}[section]
\newtheorem{lemma}[theorem]{Lemma}
\newtheorem{corollary}[theorem]{Corollary}
\newtheorem{proposition}[theorem]{Proposition}
\newtheorem{conjecture}[theorem]{Conjecture}

\theoremstyle{definition}
\newtheorem{definition}[theorem]{Definition}
\newtheorem{example}[theorem]{Example}

\newtheorem{situation}[theorem]{Situation}

\theoremstyle{remark}
\newtheorem{remark}[theorem]{Remark}

\usepackage{amsmath,amsfonts,latexsym}

\newcommand{\R}{\mathbb{R}}
\newcommand{\reals}{\mathbb{R}}

\newcommand{\complexs}{\mathbb{C}}

\newcommand{\naturals}{\mathbb{N}}

\newcommand{\integers}{\mathbb{Z}}

\newcommand{\rationals}{\mathbb{Q}}

\DeclareMathOperator{\id}{id}

\newcommand{\boundedops}{\mathcal{B}}

\newcommand{\abs}[1]{\left\lvert#1\right\rvert} 
\newcommand{\norm}[1]{\left\lVert#1\right\rVert}

\newcommand{\tensor}{\otimes}
\newcommand{\into}{\hookrightarrow}

\newcommand{\subgroup}{<}

\DeclareMathOperator{\supp}{supp}   
\DeclareMathOperator{\im}{im}      

\DeclareMathOperator{\diag}{diag}
\DeclareMathOperator{\spec}{Spec}
\DeclareMathOperator{\specpt}{Spec_{\text{\rm point}}}

\DeclareMathOperator{\tr}{tr}
\DeclareMathOperator{\pr}{pr}

\newcommand{\forget}[1]{}

\newcommand{\innerprod}[1]{\langle #1 \rangle}

{\catcode`@=11\global\let\c@equation=\c@theorem}

\allowdisplaybreaks[2]
\newcommand{\LinnelsC}{\mathcal{C}}

\newcommand{\extendedC}{{{\mathcal{D}}}}

\DeclareMathOperator{\Det}{det}
\newcommand{\NeumannN}{\mathcal{N}}
\newcommand{\universalU}{\mathcal{U}}
\newcommand{\RAgroups}{\mathcal{G}}

\renewcommand{\det}{\Det}
\begin{document}
\date{Last compiled \today; last edited  \heuteIst{} or later}

\title{Approximating $L^2$-invariants, and the Atiyah conjecture}
\author{J{\'o}zef Dodziuk\thanks{email: jdodziuk@gc.cuny.edu\protect\\
Ph.D.\ Program in Mathematics, Graduate Center of CUNY, New York, NY 10016, USA.}
\and Peter Linnell\thanks{email: linnell@math.vt.edu\protect
\\Department of Mathematics, Virginia Tech Blacksburg, VA 24061-0123, USA.} 
\and Varghese Mathai\thanks{e-mail: vmathai@math.mit.edu, vmathai@maths.adelaide.edu.au\protect\\
Department of Mathematics, University of Adelaide, Adelaide 5005, Australia,\protect\\ 
Department of Mathematics, Massachusetts Institute of Technology, Cambridge, MA 02139, USA.\protect\\
Research funded by the Clay Mathematical Institute.} 
\and Thomas Schick\thanks{
e-mail: {schick@uni-math.gwdg.de}\protect\\
FB Mathematik, Georg-August-Universit{\"a}t G{\"o}ttingen, Bunsenstr.~3, 37073 G{\"o}ttingen, Germany.\protect\\
Research in part funded by DAAD (German Academic Exchange Agency).}  
\and Stuart Yates\thanks{e-mail: syates@maths.adelaide.edu.au\protect\\
Department of Mathematics, University of Adelaide, Adelaide 5005, Australia,\protect\\ 
Research funded by the Australian Research Council.}\\
}
        
\maketitle

\begin{abstract}
Let $G$ be a torsion free discrete group and let
$\overline\rationals$ denote the
field of algebraic numbers in $\complexs$.  We prove that
 $\overline\rationals G$ fulfills  the Atiyah conjecture if $G$ lies in
 a certain class of groups $\extendedC$, which contains in
 particular all groups which are
residually
 torsion free elementary amenable or which are residually free. 
This result implies that there are no non-trivial zero-divisors
 in $\complexs G$.
The statement relies on new approximation results for $L^2$-Betti numbers over
 $\overline \rationals G$, which are the core of the work done in this 
 paper.
 Another set of results in the paper is concerned with certain number 
 theoretic properties of eigenvalues for the combinatorial Laplacian 
 on $L^2$-cochains on any normal covering space of a finite $CW$ complex. 
 We establish the absence of  eigenvalues that are transcendental numbers, 
 whenever the covering transformation group is either amenable or in the 
 Linnell class $\mathcal C$. We also establish the absence of eigenvalues 
 that are Liouville transcendental numbers whenever the covering transformation 
 group is either residually finite or more generally in
a certain large bootstrap class $\mathcal G$.

Please take the errata to Schick: ``$L^2$-determinant class and approximation
of 
$L^2$-Betti numbers'' into account, which are added at the end of the file,
rectifying some unproved statements about ``amenable 
extension''. As a consequence, throughout, amenable extensions should be
extensions with \emph{normal} subgroups.

(non-commutative measure theory)
\end{abstract}

\newpage
\section{Introduction}

Atiyah \cite{Atiyah(1976)} defined in 1976 the $L^2$-Betti numbers of
a compact
Riemannian manifold. They are defined in terms of the spectrum of the
Laplace operator on the universal covering of $M$. By Atiyah's
$L^2$-index theorem \cite{Atiyah(1976)}, they can be used e.g.~to
compute the Euler
characteristic of $M$.

Dodziuk \cite{Dodziuk(1977)} proved an $L^2$-Hodge de Rham theorem
which gives a
combinatorial interpretation of the $L^2$-Betti numbers, now in terms
of the spectrum of the combinatorial Laplacian. This operator can be
considered to be a matrix over the integral group ring of the
fundamental group $G$ of $M$. These $L^2$-invariants have been
successfully used to give new results in differential geometry,
topology, and algebra. For example, one can prove certain cases of the 
Hopf conjecture about the sign of the Euler characteristic of
negatively curved manifolds \cite{Gromov(1991)}, or, in a completely
different direction, certain cases of the zero divisor conjecture,
which asserts that $\rationals[G]$ has no non-trivial zero divisors if 
$G$ is torsion-free \cite[Lemma 2.4]{Lueck(2001)}.

In this paper, we will be concerned with the computation of the
spectrum of matrices not only over $\integers G$, but also over
$\overline\rationals G$, where $\overline\rationals$ is the field of
algebraic numbers in $\complexs$. This amends our understanding of
the combinatorial version of $L^2$-invariants. Moreover, it allows
important generalizations of the algebraic applications. In
particular, we can prove in new cases that $\complexs G$, and not
only $\rationals G$, has no non-trivial zero divisors.

Specifically, the goal of this note is to extend a number of 
(approximation) results
of \cite{Lueck(1994c)}, \cite{Dodziuk-Mathai(1998)}, 
\cite{Farber(1998)},
\cite{Schick(1998a)} and \cite{Schick(1999)} from $\rationals G$
to $\overline\rationals G$, and to give more
precise information about the spectra of matrices over
$\overline\rationals G$ than was given there. 

We deal with the following situation: $G$ is a discrete group and $K$
is a subfield of $\complexs$. The most important examples are
$K=\rationals$, $K=\overline\rationals$, and $K=\complexs$. Assume
$A\in M(d\times d, KG)$. We use the notation 
$l^2(G)=\{f\colon G\to\complexs\mid\sum_{g\in G} \abs{f(g)}^2<\infty\}$.

By left convolution, $A$ induces a bounded linear operator $A\colon
l^2(G)^d\to l^2(G)^d$ (using the canonical left $G$-action on
$l^2(G)$), which commutes with the right $G$-action.

 Let $\pr_{\ker A}\colon l^2(G)^d\to l^2(G)^d$ be the orthogonal
 projection onto $\ker A$. Then
\begin{equation*}
\dim_G(\ker A): =\tr_G(\pr_{\ker A}):=\sum_{i=1}^n \innerprod{\pr_{\ker A}
  e_i,e_i}_{l^2(G)^n},
\end{equation*}
where $e_i\in l^2(G)^n$ is the vector with the trivial element of
$G\subset l^2(G)$
at the $i^{\text{th}}$-position and zeros elsewhere.

To understand the type of results we want to generalize, we repeat the
following definition:

\begin{definition}\label{def:AtiyahCon} Let $K$ be a subring of $\complexs$.
  We say that a torsion free group $G$ fulfills the
 \emph{strong Atiyah
conjecture over $KG$} if 
\begin{equation*}
              \dim_G (\ker A)\in \integers \qquad\forall A\in
M(d\times d, K G);
              \end{equation*}
where $\ker A$ is the kernel of the induced map $A\colon l^2(G)^d\to
l^2(G)^d$.
\end{definition}
It has long been observed that for a torsion free group $G$, 
the strong Atiyah conjecture over
$K G$ implies that $K G$ does not contain non-trivial
zero divisors. If $K=\overline \rationals$, it does even imply that
$\complexs G$ has no non-trivial zero divisors (we will recall the
argument for this in Section \ref{sec:algebraic_to_complex}).

Linnell proved the strong Atiyah conjecture over $\complexs G$ if $G$
is torsion free and
$G\in\LinnelsC$, where $\LinnelsC$ denotes the smallest class of groups
containing all free groups and which is
closed under extensions with elementary amenable quotient and under
directed unions.

In \cite{Schick(1999)} (see also \cite{Schick(2001)}), Linnell's
results have been generalized to a larger class $\extendedC$ of
groups (see Definition \ref{def:extC} below), but only for $\rationals G$
instead of $\complexs G$. One of
the objectives of this paper is, to fix the flaw that only the
coefficient ring $\rationals$ is allowed.

Recall that the class of elementary amenable groups is the smallest
class of groups containing all cyclic and all finite groups and which
is closed under taking group extensions and directed unions. 
The class $\extendedC$ as defined in \cite{Schick(1999)}
is then as follows.


\begin{definition}\label{def:extC}
  Let $\extendedC$ be the smallest non-empty class of groups such that:
  \begin{enumerate}
  \item \label{aex} If $G$ is torsion-free, $A$ is elementary
    amenable, and we have an epimorphism
    $p\colon G\to A$ such that $p^{-1}(E)\in\extendedC$ for every
    finite subgroup $E$ of $A$, then $G\in\extendedC$.
   \item \label{sgr} $\extendedC$ is subgroup closed.
   \item\label{lim} Let $G_i\in \extendedC$ be a
  directed system of groups
    and $G$ its (direct or inverse) limit. Then $G\in\extendedC$.
\end{enumerate}
\end{definition}

Certainly $\extendedC$ contains all free groups, because these are
residually torsion-free nilpotent.
One of the main results of \cite{Schick(1999)} (see also
\cite{Schick(2001)}) is
\begin{theorem}\label{theo:strongAt}
  If  $G\in\extendedC$ then the strong Atiyah
  conjecture over 
  $\rationals G$ is true.
\end{theorem}

We improve this to
\begin{theorem}\label{theo:algebraicstrongAt}
  If $G\in\extendedC$ then the strong Atiyah conjecture over
  $\overline\rationals G$ is true.
\end{theorem}

As mentioned above, as a corollary we obtain via Proposition
\ref{prop:algebraic_to_complex} that there are no non-trivial zero
divisors in $\complexs G$, if $G\in\extendedC$.

\begin{corollary}\label{corol:AC_res_tf_nilpotent}
  Suppose $G$ is residually torsion free nilpotent, or residually
  torsion free solvable. Then the strong Atiyah conjecture is true for
  $\overline\rationals G$.
\end{corollary}
\begin{proof}
  Under the assumptions we make, $G$ belongs to $\extendedC$.
\end{proof}

The proof of Theorem \ref{theo:strongAt} relies on certain
approximation results proved in \cite{Schick(1998a)} which give
information about the spectrum of a self-adjoint matrix $A\in
M(d\times d,\rationals G)$. We have to improve these results to
matrices over $\overline\rationals G$. This will be done in Section
\ref{sec:algapr},
where also the precise statements are given. As a special case, we
will prove the following result.
\begin{theorem}\label{theo:residual_approxi_algebraic_coeff}
  Suppose the group $G$ has a sequence of normal subgroups $G\supset G_1\supset
  G_2\supset \cdots$ with $\bigcap_{i\in\naturals} G_i=\{1\}$ and such that 
 $G/G_i\in\extendedC$ for every $i\in\naturals$. Assume that $B\in
 M(d\times d,\overline\rationals G)$. Let $p_i\colon G\to G/G_i$ denote
 the natural epimorphism, and let $B[i]\in M(d\times d,\overline\rationals
 G/G_i)$ be the image of $B$ under the ring homomorphism induced by
 $p_i$. Then
 \begin{equation*}
   \dim_G(\ker(B)) = \lim_{i\to\infty} \dim_{G/G_i} \ker(B[i]).
 \end{equation*}
\end{theorem}

We will also have to prove appropriate approximation results for
generalized amenable extensions, whereas in \cite{Schick(1998a)} we
only deal with the
situation $H\subgroup G$ such that $G/H$ is an amenable homogeneous
space. Actually, in the amenable situation we are able to generalize
F{\o}lner type approximation results to all matrices over the complex
group ring.

Consider a larger class $\RAgroups$ of groups defined as follows,
following the definitions in \cite{Schick(1999)}.

\begin{definition}\label{def:generalized_amenable_extension}
  Assume that $G$ is a finitely generated discrete group with a finite
  symmetric set of generators $S$ (i.e.~$s\in S$ implies $s^{-1}\in
  S$), and let $H$ be an
  arbitrary discrete group. We say that \emph{$G$ is a generalized
    amenable extension of $H$}, if there is a set $X$ with a free $G$-action
  (from the left) and a commuting free $H$-action (from the right),
  such that a sequence of $H$-subsets $X_1\subset X_2\subset
  X_3\subset \cdots \subset X$ exists with $\bigcup_{k\in\naturals} X_k =X$,
  and with $\abs{X_k/H}<\infty$ for every $k\in\naturals$, and such
  that
  \begin{equation*}
   \frac{ \abs{(S\cdot X_k - X_k)/H}}{\abs{X_k/H}}
   \xrightarrow{k\to\infty} 0.
 \end{equation*}
\end{definition}

\begin{definition}
  Let $\RAgroups$ be the smallest class of groups which contains the trivial
  group and is closed under the following processes:
  \begin{enumerate}
   \item If $H\in\RAgroups$ and $G$ is a generalized amenable
     extension of $H$, then $G\in\RAgroups$.
  \item If $H\in \RAgroups$ and $U\subgroup H$, then $U\in\RAgroups$.
  \item If $G=\lim_{i\in I}G_i$ is the direct or inverse limit of a
    directed system of groups $G_i\in \RAgroups$,
    then  $G\in\RAgroups$.
  \end{enumerate}
\end{definition}
We have the inclusion $\LinnelsC\subset\extendedC\subset\RAgroups$,
and in particular the class $\RAgroups$ contains all amenable groups,
free groups, residually finite groups, and residually amenable groups. 

The approximation results derived in sections \ref{sec:algapr} and
\ref{sec:complapr} then imply the following ``stability'' result about the
Atiyah conjecture.
\begin{proposition}
  \label{prop:algebraiclimits}
  Assume $G$ is a subgroup of the direct or inverse limit of a
  directed system of groups $(G_i)_{i\in I}$. Assume
  $G_i\in\RAgroups$, $G_i$ is torsion-free, and $G_i$ satisfies the
  strong Atiyah conjecture over $\overline\rationals G_i$ for every
  $i\in I$.  Then $G$ also satisfies the strong Atiyah conjecture over
  $\overline\rationals G$.
\end{proposition}

In Section
\ref{sec:complapr} we will address the question when Theorem
\ref{theo:residual_approxi_algebraic_coeff} holds for all matrices
over $\complexs G$. We will answer this
question affirmatively in particular if $G$ is torsion free and
elementary amenable.

\begin{proposition}
  \label{prop:resapproxi}
  Let $G$ be a torsion free elementary amenable group with a nested
  sequence of normal subgroups $G\supset G_1\supset G_2\supset\cdots$
  such that $\bigcap G_k=\{1\}$. Assume $G/G_i\in\RAgroups$ for every 
  $i\in\naturals$. Let $A\in M(d\times d,\complexs G)$
  and $A_k\in M(d\times d,\complexs [G/G_k])$ its image under the maps
  induced from the epimorphism $G\to G/G_k$. Then
  \begin{equation*}
    \dim_G(\ker(A)) = \lim_{k\to\infty} \dim_{G/G_k}(\ker(A_k)).
  \end{equation*}
\end{proposition}

In section 6 we study the algebraic eigenvalue property for discrete 
groups $G$. This implies certain number theoretic properties of eigenvalues
for  operators $A = B^*B$, where \hbox{$B \in M(d\times d,\overline\rationals G)$};
this includes the combinatorial Laplacian on $L^2$-cochains on any normal covering
space of a finite $CW$ complex. We establish the absence of  eigenvalues
that are transcendental numbers, whenever the 
group is either amenable or in the Linnell class $\LinnelsC$. We also
establish the absence of eigenvalues that are Liouville transcendental
numbers whenever the group is either residually
finite or more generally in the bootstrap class $\RAgroups$. The 
statements rely on new approximation results for spectral density
functions of such operators $A$.

\section{Notation and Preliminaries}

Let $G$ be a (finitely generated) group and $A\in M(d\times d,\complexs G)$.
We consider $A$ as a matrix indexed
by $I\times I$ with $I=\{ 1,\ldots,d\}\times G$. We can equally
well regard $A$ as a family of $d \times d$ complex matrices indexed by
$G \times G$. With this interpretation, the following
invariance condition holds 
\begin{equation} \label{eq:inv-cond}
A_{x,y} = A_{xz,yz}.
\end{equation}
for all $x,y,z \in G$.
Similarly, elements of 
$l^2(G)^d$ will be regarded as sequences of elements of 
${\mathbb C}^d$ indexed by elements of $G$ so that 
\begin{equation*}
(Af)_x = \sum_y A_{x,y} f_y.
\end{equation*}

Note that one can also identify $l^2(G)^d$ with
$l^2(G)\tensor_\complexs \complexs^d$, and similarly for matrices, as
we will do occasionally.

Observe also that there exists a positive number $r$ depending on $A$
only such that 
\begin{equation}\label{eq:sparse-matrix}
A_{x,y}=0 \qquad\qquad\mbox{if}\qquad \rho(x,y)>r
\end{equation}
where $\rho(x,y)$ denotes the distance in the word metric between 
$x$ and $y$ for a given finite set of generators of $G$.

\subsection{Free actions}
\label{sec:free_actions}

If $G$ acts freely on a space $X$ (as in Definition
\ref{def:generalized_amenable_extension}), then every matrix $A$ as
above (and actually every operator on $l^2(G)^d$) induces an operator
$A_X$ on $l^2(X)^d$. 

Choosing a free $G$-basis for $X$, $l^2(X)^d$ can be identified with a
(possibly infinite) direct sum of copies of $l^2(G)^d$, and $A_X$ is a
diagonal operator with respect to this decomposition, where each
operator on $l^2(G)^d$ is equal to $A$. Consequently, if we apply the
functional calculus, $f(A)_X=f(A_X)$ for every measurable function $f$ 
(provided $A$ is self-adjoint).

Pick $x\in X$ and define 
\begin{equation*}
\tr_G(A_X):= \sum_{i=1}^d \innerprod{A_X (x\tensor
  e_i),x\tensor e_i}=\sum_{i=1}^d \innerprod{A (1\tensor e_i),1\tensor e_i},
\end{equation*}
where $(e_1,\dots,e_d)$ is the standard basis
of $\complexs ^d$ and where we identify $l^2(X)^d$ with a direct sum of 
a number of copies of $l^2(G)^d$ as above and use the identification
$l^2(X)^d=l^2(X)\tensor_\complexs \complexs^d$. Then
$\tr_G(A_X)=\tr_G(A)$ (in
particular the expression does not depend on $x\in X$).

\subsection{Spectral density functions and determinants}
\label{sec:spectral_density}

Assume $A\in M(d\times d,\complexs G)$ as above, and assume that $A$,
considered as an operator on $l^2(G)^d$, is positive self-adjoint.
This is the case e.g.~if $A=B^*B$ for some $B\in M(d\times
d,\complexs G)$. 
Here $B^*$ is the adjoint in the sense of operators on $l^2(G)^d$. If
$B=(b_{ij})$ with $b_{ij}\in \complexs G$, then
 $B^*\in M(d\times d,\complexs G)$ is the matrix with entry $b_{ij}^*$ 
 at the position $(j,i)$, where $(b_{ij}^*)_{g}=
 \overline{b_{ij}}_{g^{-1}}$.

For $\lambda\ge 0$ let $\pr_\lambda$ be the spectral projection of $A$
corresponding to the interval $[0,\lambda]$. This operator also
commutes with the right $G$-action on $l^2(G)^d$. Using exactly the
same definition of the regularized trace as before for finite
matrices, we define the spectral density function
\begin{equation*}
  F_A(\lambda):= \tr_G(\pr_\lambda).
\end{equation*}
This function is right continuous monotonic increasing, with
\begin{equation*}
F_A(0)=\dim_G(\ker(A)).
\end{equation*}
We will be interested in the behavior of
this function near $0$.

Using the spectral density function, we can define a normalized
determinant as follows.
\begin{definition}\label{def:det_G}
  If $A\in M(d\times d,\complexs G)$ is positive self-adjoint, we
  define (the logarithm) of its \emph{normalized determinant} by the
  Riemann-Stieltjes integral
  \begin{equation*}
    \ln\det_G(A) = \int_{0^+}^\infty \ln(\lambda) \; dF_A(\lambda).
  \end{equation*}
  Using integration by parts, if $K\ge \norm{A}$, we can rewrite
  \begin{equation*}
    \ln\det_G(A) = \ln(K)(F_A(K)-F_A(0))-\int_{0^+}^K
    \frac{F_A(\lambda)-F_A(0)}{\lambda} \;d\lambda.
  \end{equation*}
\end{definition}
Observe that, if $G=\{1\}$, $\det_G(A)$ is the product of the
eigenvalues of $A$ which are different from zero. In general,
$\det_G(A)$ is computed from the spectrum of $A$ on the orthogonal
complement of its kernel.

\subsection{Generalized amenable extensions}

Let $G$ be a generalized amenable extension of $H$ with
a finite set of symmetric
generators $S$. Let $X_1\subset X_2\subset \cdots
\subset X$ be a generalized
amenable exhaustion as in Definition
\ref{def:generalized_amenable_extension}. 

\begin{definition}
   For $r\in\naturals$, we define the $r$-neighborhood of the boundary
 of $X_i$ (with respect to $S$) as
 \begin{equation*}
   \mathfrak{N}_r(X_i):= S^r X_i\cap S^r(X-X_i).
 \end{equation*}
 Note that all these sets are right $H$-invariant.
\end{definition}

\begin{lemma}\label{lem:amenable_exhaustion}
  We have
  \begin{align}
     & \frac{ \abs{(S^n\cdot X_k - X_k)/H}}{\abs{X_k/H}}
     \xrightarrow{k\to\infty} 0\qquad\forall n\in\naturals, \label{eq:1}\\
     & \frac{ \abs{(\mathfrak{N}_n(X_k)/H)}}{\abs{X_k/H}}
     \xrightarrow{k\to\infty} 0\qquad\forall n\in\naturals.\label{eq:4}
\end{align}
 In particular, for any finite subset $T$ of $G$,
 \begin{equation}\label{eq:2}
      \frac{ \abs{(T\cdot X_k \cap T(X- X_k))/H}}{\abs{X_k/H}}
   \xrightarrow{k\to\infty} 0.
 \end{equation}
\end{lemma}
\begin{proof}
  Observe that $S^n X_k = S^{n-1} X_k\cup S(S^{n-1} X_k -
  X_k)$. Consequently,
  \begin{equation*}
    \frac{\abs{(S^n X_k-X_k)/H}}{\abs{X_k/H}} \le \frac{\abs{(S^{n-1}
          X_k-X_k)/H}}{\abs{X_k/H}} + \abs{S}
      \frac{\abs{(S^{n-1}X_k-X_k)/H}}{\abs{X_k/H}}.
    \end{equation*}
    Since $\abs{S}$ is fixed, \eqref{eq:1} follows by
    induction. 

  Next, $\mathfrak{N}_n(X_k)=S^n X_k\cap S^n(X-X_k)$ is contained in the union of $S^n
  X_k-X_k$ and $S^n (X-X_k)\cap X_k$. For the latter observe that
  $vy=x$ with $v\in S^n$, $y\notin X_k$ and $x\in X_k$ means that
  $y=v^{-1}x$, and since $S$ is symmetric $v^{-1}\in
  S^n$, i.e.~$y\in S^nX_k-X_k$. Consequently, $\abs{(S^n (X-X_k)\cap
    X_k)/H}\le
  \abs{S^n}\cdot \abs{(S^n X_k-X_k)/H}$, such that \eqref{eq:4}
    follows from \eqref{eq:1}.

\eqref{eq:2} follows, since $T\subset S^n$ for
    sufficiently large $n\in\naturals$.
\end{proof}

\begin{example}\label{ex:amanble_action}
  If $H\subset G$ and $G/H$ is an amenable Schreier graph (e.g.\ if $H$
  is a normal subgroup of $G$ such that $G/H$ is amenable), 
  then $G$ is a generalized amenable extension of $H$. In
  these cases, we can take $X$ to be $G$, and $X_k:= p^{-1} Y_k$ if
  $p\colon G\to G/H$ is the natural epimorphism, and $Y_k$ are a
  F{\o}lner exhaustion of the  amenable Schreier graph (or amenable
  group,
  respectively) $G/H$.
\end{example}

The following result follows from an idea of Warren Dicks.
\begin{example}
  Assume $G$ is a finitely generated group with symmetric set of
  generators $S$ which acts on a set $X$ (not
  necessarily freely). Assume there is an exhaustion $X_1\subset
  X_2\subset \cdots \subset X$ of $X$ by finite subsets such that 
  \begin{equation*}
    \frac{\abs{S X_k-X_k}}{\abs{X_k}}\xrightarrow{k\to\infty} 0.
  \end{equation*}
  For $x\in X$ let $G_x:=\{g\in G\mid gx=x\}$ be the stabilizer of
  $x$. If a group $H$ contains an abstractly isomorphic copy of $G_x$
  for every $x\in X$, e.g.~if $H$ is the free product or the direct
  sum of all $G_x$ ($x\in X$), then $G$ is a generalized amenable
  extension of $H$.
\end{example}

This follows immediately from the following Lemma, which is (with its 
proof) due to Warren Dicks.

\begin{lemma}\label{lem:action}
   Let $X$ be a left $G$-set.

   If, for every $x \in X$, the $G$-stabilizer $G_x$ embeds in $H$,
   then
\begin{enumerate}
\item\label{item:action} there exists a family of maps
   $$(\gamma_x\colon G \to H \mid x \in X)$$
   such that, for all $x \in X$, for all $g_1$, $g_2 \in G$,
   $$\gamma_{g_1x}(g_2)\gamma_{x}(g_1)= \gamma_{x}(g_2g_1),$$
   and the restriction of $\gamma_x$ to $G_x$ is injective
   (and a group homomorphism).
\end{enumerate}
   If~\ref{item:action} holds, then there is a $G$-free $H$-free
   $(G,H)$-bi-set structure on $X \times H$ such that
   $$g(x,y)h = (gx, \gamma_x(g)yh),$$ for all
   $g \in G$, $x \in X$,  $y,h \in H$.
   Moreover, for each $x \in X$, if $Y$ is a left
   $\gamma_x(G_x)$-transversal in $H$, then $\{x\}\times Y$ is a
   left $G$-transversal in $Gx \times H$.
\end{lemma}

\begin{proof}
   Let $x_0 \in X$, and let $\alpha\colon G_{x_0} \to H$ be an injective
   group homomorphism.

   We claim that there exists a right $\alpha$-compatible map
   $\beta\colon G \to H$,  that is, $\beta(gk) = \beta(g)\alpha(k)$ for
   all $k \in G_{x_0}$ and all $g \in G$.  Since $G$ is free as right
   $G_{x_0}$-set, we can construct $\beta$ as follows.  Choose a right
   $G_{x_0}$-transversal in $G$, and map each element of this transversal
   arbitrarily to an element of $H$, and extend by the $G_{x_0}$-action
   to all of $G$.

   Now let $x \in Gx_0$, and choose $g_x \in G$ such that $x = g_x x_0$.

   Define $\gamma_x\colon G \to H$ by
   $\gamma_x(g) = \beta(gg_x)\beta(g_x)^{-1}$ for all $g\in G$.
   Since $\beta$ is right $\alpha$-compatible, we see that $\gamma_x$ is
   independent of the choice of $g_x$.

   For all $g_1$, $g_2 \in G$,
   \begin{align*}
   \gamma_{g_1x}(g_2)\gamma_{x}(g_1)&=
   \beta(g_2g_1g_x)\beta(g_1g_x)^{-1}\beta(g_1g_x)\beta(g_x)^{-1}\\
   &= \beta(g_2g_1g_x)\beta(g_x)^{-1} = \gamma_{x}(g_2g_1)
   \end{align*}

   We now show that $\gamma_x$ restricted to $G_x$ is an injective group
   homomorphism.   Suppose that $g \in G_x$, so
   $g_x^{-1}gg_x \in G_{x_0}$, and
  \begin{align*}
   \gamma_x(g) &= \beta(gg_x)\beta(g_x)^{-1}
    = \beta(g_x g_x^{-1}gg_x)\beta(g_x)^{-1} \\
    &= \beta(g_x)\alpha(g_x^{-1}gg_x)\beta(g_x)^{-1}.
   \end{align*}
   Since $\alpha$ is injective, we see that $\gamma_x$ is injective on
   $G_x$.

   Since $Gx_0$ is an arbitrary $G$-orbit in $X$, we have proved
   that~\ref{item:action} holds.

   Now suppose that~\ref{item:action} holds.  It is straightforward
   to check that we have the desired $(G,H)$-bi-set structure on
   $X\times H$, and that it is $G$-free and $H$-free.  Finally, for any
   $x \in X$, there is a well-defined, injective map from
   $\gamma_x(G_x)\backslash H$ to $G\backslash(X\times H)$,  with
   $\gamma_x(G_x)h$ mapping to $G(x,h)$, for all $h \in H$. It is
   surjective, since
   $G(gx,h) = G(x,\gamma_x(g)^{-1}h)$ for all $g \in G$.
\end{proof}

One should observe that such a $G$-space $X$ of course is the disjoint 
union of Schreier graphs $G/G_x$. However, $X$ being an amenable
$G$-set in the sense of Example \ref{ex:amanble_action} does not
necessarily imply that any of the these $G/G_x$ is an amenable
Schreier graph. If, for a fixed set of generators $S$ of $G$ one can
find quotient groups with arbitrarily small exponential growth, then
the disjoint union of these quotients becomes an amenable $G$-set.

\subsection{Direct and inverse limits}\label{limits}

 Suppose that a group $G$ is the direct or inverse limit of a
  directed system of groups $G_i$, $i\in I$. The latter means that we have
  a partial ordering $<$ on $I$ such that for all $i,j\in I$ there is $k\in
  I$ with $i<k$ and $j<k$, and maps $p_{ij}\colon G_i \to G_j$ in the
case of the direct limit and $p_{ji} \colon G_j \to G_i$ in the case
of the inverse limit whenever $i<j$, 
  satisfying the obvious compatibility conditions. In the case of a 
  direct limit we let
  $p_i\colon G_i\to G$ be the natural maps and similarly,  for
  an inverse limit, $p_i\colon G\to G_i$.
The maps between groups induce
  mappings between matrices with coefficients in group rings in the
  obvious way. Let  $B\in M(d\times d,\complexs G)$. If $G$ is an 
  inverse limit
  we set $B[i]=p_i(B)$. 
  
  In the case of a direct limit it is necessary to make
  some choices. Namely, let
$B=(a_{kl})$ with $a_{kl}=\sum_{g\in G}\lambda^g_{kl}g$. Then, only
finitely many of the $\lambda^g_{kl}$ are nonzero. Let $V$ be the
corresponding finite collection of $g\in G$. Since $G$ is the
direct limit of $G_i$, we can find $j_0\in I$ such that $V\subset
p_{j_0}(G_{j_0})$. Choose an inverse image for each $g$ in
$G_{j_0}$. This gives a matrix $B[j_0]\in M(d\times
d,G_{j_0})$ which is mapped to 
$B[i]:=p_{j_0i}B[j_0]\in M(d\times d,G_i)$ for
$i>j_0$. 

Observe that in both cases $B^*[i]=B[i]^*$.

\section{Approximation with algebraic coefficients}
\label{sec:algapr}

\subsection{Lower bounds for determinants}
\label{sec:determinant_bounds}

We will define in Definition \ref{def:kappa} an invariant $\kappa(A)$
of matrices $A\in M(d\times d,\complexs G)$ which occurs in the lower
bounds of generalized determinants we want to establish.

\begin{definition}
  \label{def:detestprop}
  We say that a group $G$ has the \emph{determinant bound property} if,
  for every $B\in M(d\times d,o(\overline\rationals)
  G)$ and $\Delta=B^*B$, where $o(\overline\rationals)$ is the ring of
  algebraic 
  integers over $\integers$ in $\complexs$, the following is true:
  Choose a finite Galois
  extension $i\colon L\subset \complexs$ of $\rationals$ such that $B\in
  M(d\times d,LG)$. Let $\sigma_1,\dots,\sigma_r\colon L\to\complexs$ be the
  different embeddings of $L$ in $\complexs$ with $\sigma_1=i$. Then
  \begin{equation}\label{basicest}
    \ln\det_G(\Delta) \ge -d \sum_{k=2}^r
    \ln(\kappa(\sigma_k(\Delta))).
  \end{equation}

  We say that $G$ has the \emph{algebraic continuity property}, if always
 \begin{equation}\label{sigmaequal}
    \dim_G(\ker(B))= \dim_G(\ker(\sigma_k B))\quad\forall k=1,\dots,r.
  \end{equation}
\end{definition}
\begin{theorem}\label{newdetest}
   If  $G\in\RAgroups$, then $G$ has the determinant bound property
   and the algebraic continuity property.
\end{theorem}

\begin{corollary}
  If $B=B^*\in M(d\times d,\rationals G)$ and $\lambda\in\reals$ is
  algebraic and an eigenvalue of $B$, then $\lambda$ is totally real,
  i.e.~$\sigma(\lambda)\in\reals$ for every automorphism
  $\sigma\colon\complexs\to \complexs$. The same is true if
  $\rationals$ is replaced by any totally real algebraic extension
  $\rationals\subset R\subset \complexs$,
  i.e.~$\sigma(R)\subset\reals$ for every automorphism
  $\sigma\colon\complexs\to\complexs$. 
\end{corollary}

Estimate \eqref{basicest} implies a sometimes more convenient estimate 
for the spectral density functions, which we want to note now.
\begin{corollary}
Let $C \in \mathbb {R}$.
  Suppose $A\in M(d\times d,\complexs G)$ with $\norm{A}\ge 1$ is
  positive self-adjoint and
  satisfies 
  \begin{equation*}
    \ln\det_G(A)\ge -C.
  \end{equation*}
  Then, for $0<\lambda<\norm{A}$,
  \begin{equation}\label{eq:6}
    F_A(\lambda)-F_A(0)\le
    \frac{C+d\ln(\norm{A})}{-\ln(\lambda/\norm{A})}. 
  \end{equation}
  In particular, in the situation of Theorem \ref{newdetest},
  for all $\lambda$ such that $0<\lambda<\norm{A}$, we have
  \begin{equation}\label{eq:7}
    \begin{split}
           F_A(\lambda)-F_A(0) & \le \frac{d\cdot (\ln(\norm{A})
             +\sum_{k=2}^r
             \ln(\kappa(\sigma_k(A))))}{-\ln(\lambda/\norm{A})}\\
           & \le
      \frac{d\sum_{k=1}^r \ln(\kappa(\sigma_k(A)))}{-\ln(\lambda/\norm{A})}.
  \end{split}
\end{equation}
If $r=1$, i.e.~$A\in M(d\times d,\integers G)$, then this  simplifies
to 
\begin{equation}\label{eq:8}
  F_A(\lambda)-F_A(0) \le
  \frac{d\ln(\norm{A})}{-\ln(\lambda/\norm{A})}.
\end{equation}
\end{corollary}
\begin{proof}
  The proof of \eqref{eq:6} is done by an elementary estimate of
  integrals. Fix $0<\lambda<\norm{A}$. Then
  \begin{equation*}
    \begin{split}
     -C & \le  \ln\det_G(A) = \int_{0^+}^{\norm{A}} \ln(\tau)
     \;dF_A(\tau)\\
      &= \int_{0^+}^\lambda \ln(\tau) \; dF_A(\tau) +
      \int_{\lambda^+}^{\norm{A}} \ln(\tau)\; dF_A(\tau)\\
      &\le \ln(\lambda) (F_A(\lambda)-F_A(0)) + \ln(\norm{A})(
      F_A(\norm{A})-  F_A(\lambda))\\
      & \le 
      \ln(\lambda/\norm{A}) (F_A(\lambda)-F_A(0)) + d \ln(\norm{A})
     .
    \end{split}
  \end{equation*}
  This implies that $$
   F_A(\lambda) -F_A(0)  \le \frac{C +
        d\ln(\norm{A})}{-\ln(\lambda/\norm{A}) }.$$
        
  For \eqref{eq:7} we use \eqref{basicest} and the fact that
  $\norm{A}\le \kappa(\sigma_1(A))$.
\end{proof}

\subsection{Approximation results}\label{sec:appr-results}

We will prove Theorem \ref{newdetest} together with certain approximation
results for $L^2$-Betti numbers in an inductive way. We now 
formulate precisely the approximation results, which were already
mentioned in the introduction. To do
this, we first describe the situation that we will consider.

\begin{situation}\label{newsit}
  Given is a group $G$ and a matrix $B\in M(d\times d,\complexs
  G)$. We assume that one of the following applies:
  \begin{enumerate}
  \item $G$ is the inverse limit of a directed system of groups $G_i$. Let
    $B[i]\in M(d\times d,\complexs G)$
    be the image of $B$ under the natural map $p_i: G\to G_i$.
  \item $G$ is the direct limit of a directed system of groups
    $G_i$. For $i\ge i_0$ we choose $B[i]\in M(d\times d,\complexs G)$
    as described in subsection \ref{limits}.
  \item $G$ is a generalized amenable extension of $U$ with free
    $G$-$U$ space $X$ and F{\o}lner exhaustion $X_1\subset X_2\subset
    \cdots \subset X$ of $X$ as in Definition
    \ref{def:generalized_amenable_extension} (with
    $\abs{X_i/U}:=N_i<\infty$). Let $P_i\colon l^2(X)^d\to
    l^2(X_i)^d$ be the corresponding projection operators. Recall that $B_X$
    is the operator on
    $l^2(X)^d$ induced by $B$, using the free action of $G$ on
    $X$. We set $B[i]:= P_iB_X P^*_i$.
  \end{enumerate}

  Moreover, we set $\Delta:=B^*B$, and $\Delta[i]=B[i]^*B[i]$ (for
  $i\ge i_0$).

  In the first two cases we write $\dim_i$, $\tr_i$ and $\det_i$ for
  $\dim_{G_i}$, $\tr_{G_i}$ and $\det_{G_i}$, respectively. In the
  third case we use a slight modification. We set
  $\dim_i:=\frac{1}{N_i}\dim_U$ and $\tr_i:=\frac{1}{N_i}\tr_U$, 
  $\ln\det_i:=\frac{1}{N_i}\ln\det_{U}$, and
  write $G_i:=U$ to unify the notation.

  We also define the spectral density functions
  $F_{\Delta[i]}(\lambda)$ as in Subsection
  \ref{sec:spectral_density}, but using $\tr_i$ instead of the trace
  used there.
\end{situation}

We are now studying approximation results for the dimensions of
eigenspaces. We begin with a general
approximation result in the amenable case. In particular, we
generalize the main result of \cite{Dodziuk-Mathai(1998)} from the
rational group ring to the complex group ring of an amenable group,
with a simpler proof. Eckmann \cite{Eckmann(1999)} uses
related ideas in his proof of the main theorem of
\cite{Dodziuk-Mathai(1998)}.

\begin{theorem}\label{theo:general_amenable_approxi}
  Assume $G$ is a generalized amenable extension of $U$, and
  $B\in M(d\times d, \complexs G)$. Adopt the notation and
  situation of \ref{newsit}(3). Then
  \begin{equation*}
    \dim_G(\ker(B)) = \lim_{i\to\infty} \dim_i\ker(B[i]) =
    \lim_{i\to\infty} \frac{\dim_U\ker(B[i])}{N_i}.
  \end{equation*}
\end{theorem}

We can generalize this to \ref{newsit}(1) and \ref{newsit}(2), but,
for the time being, only at the expense of restricting the coefficient 
field to $\overline\rationals$. 

\begin{theorem}\label{newapproxi}
  Suppose $B\in M(d\times d,\overline\rationals G)$ is the limit of
  matrices $B[i]$ over $\overline\rationals G_i$ as described in
  \ref{newsit}. Assume $G_i\in\RAgroups$ for all $i$. Then
  \begin{equation*}
    \dim_G(\ker(B)) = \lim_{i\in I} \dim_i(\ker(B[i])).
  \end{equation*}
\end{theorem}
Note that it is sufficient to prove Theorem \ref{newapproxi} for
$\Delta=B^*B$, since the kernels of $B$ and of $\Delta$ coincide, and since by
Lemma \ref{ringherit} $B[i]^*=(B^*)[i]$. Moreover, every algebraic number
is the product of a rational number and an algebraic integer (compare
\cite[15.24]{Dummit-Foote(1999)}). Hence we can multiply $B$ by a suitable
integer $N$ (without changing the kernels) and therefore  assume that
$B$, and $\Delta$, is a matrix over
$o(\overline\rationals)G$.

We start with one half of the proof of Theorem
\ref{theo:general_amenable_approxi}. 

\begin{lemma}\label{lem:half_proof}
  In the situation of Theorem \ref{theo:general_amenable_approxi},
  \begin{equation*}
      \dim_G(\ker(B)) \le \liminf_{i\to\infty}
  \frac{\dim_U\ker(B[i])}{\abs{X_i/U}}. 
  \end{equation*}
\end{lemma}
\begin{proof}
First observe that, in the situation of \ref{newsit}(3),
we have to discuss the relation between $B$ and
$B_X$. Since the action of $G$ on $X$ is free, $B_X$ is defined not
only for $B\in M(d\times
d,\complexs G)$, but even for $B\in M(d\times d,\NeumannN G)$, where
$\NeumannN G$ is the group von Neumann algebra of $G$, i.e.~the weak
closure of $\complexs G$ in $\boundedops(l^2(G))$, or equivalently
(using the bicommutant theorem)
the set of all bounded operators on $l^2(G)$ which commute with the
right $G$-action.

  Let $P$ be the orthogonal projection onto $\ker(B)$. Then
  $\dim_G(\ker B)=\tr_G(P)= \tr_G(P_X)$ (where $P_X$ is the orthogonal 
  projection onto $\ker(B_X)$, as explained above).

  Let $P_i$ be the orthogonal projection of $l^2(X)^d$ onto the
set of functions with values in $\complexs^d$ and supported on $X_i$.
We need to compare $P_iB_X$ and $B_XP_i$. Let $T$ be the support of $B$,
i.e.~the set of all $g\in G$ such that the coefficient of $g$ in at
least one entry in the matrix $B$ is non-zero, and fix $r\in\naturals$ 
such that $T\subseteq S^r$.  Also let $A$ be the matrix of $B_X$ with
respect to the basis $X$.  A calculation shows that
\begin{equation}\label{diff-bdry}
(P_iB_Xf - B_X P_if)_x   =\left\{
   \begin{array}{ll}
    \sum_{y \in X\setminus X_i}A_{x,y}f_y ,   
       &  \hbox{if }x \in X_i \\[2ex]
    \sum_{y \in X_i} A_{x,y}f_y, 
       & \hbox{if } x \not\in X_i.
   \end{array}
 \right.
\end{equation} 
It follows that the value of the difference is determined by the values
$f_y$ on the $r$-neighborhood $\mathfrak{N}_r(X_i)$ of the boundary of $X_i$.

We have
$$
\tr_G  P_X = \frac{1}{|X_i/U|} 
\tr_U P_i P_X
$$
by the invariance property (\ref{eq:inv-cond}).
Since $\norm{P_iP } \leq 1$, $\tr_U P_iP  \leq 
\dim_U \overline{\im P_iP} $ and $\im P_iP=P_i(\ker(B_X))$, we have
$$
\dim_G  \ker(B_X) \leq \frac{1}{|X_i/U|} \dim_U \overline{ P_i(\ker(B_X))} .
$$
To simplify the notation we identify the functions
on $X $ supported on $Y\subset X$ with functions on
$X$. With this convention, $P_i=P_i^*$ and $B[i]=P_i B_X P_i$.

Consider $f \in \ker(B_X)$, i.e.~$B_X f=0$. It follows that
$$
B[i] P_i f   =  P_i B_X P_i f - P_i B_X f =P_i(B_XP_if - P_i B_X f).
$$
By (\ref{diff-bdry}), the right-hand side is determined uniquely by
the restriction of $f$ to $\mathfrak{N}_r(X_i)$. Thus, if the restriction of
$f\in \ker(B_X)$ to $\mathfrak{N}_r(X_i)$ is zero, then
\begin{equation}\label{eq:5}
B[i] P_i f  = 0.
\end{equation}

The bounded $U$-equivariant operator $P_i\colon l^2(X)^d\to
l^2(X_i)^d$ restricts to a
bounded operator $P_i\colon \ker(B_X)\to l^2(X_i)^d$. Let
$V\subset\ker(B_X)$ be the orthogonal complement of the kernel of this 
restriction inside $\ker(B_X)$. Then $P_i$ restricts to an injective
map $V\to l^2(X_i)$ with image $P_i(\ker(B_X))$. Consequently,
$\dim_U(V)= \dim_U \overline{P_i(\ker(B_X))}$ \cite[see \S
2]{Lueck(2001)}.

Let $\pr_i\colon l^2(X_i)^d\to \ker(B[i])$
be the orthogonal projection, and let
$Q_i\colon l^2(X)^d\to l^2(\mathfrak{N}_r(X_i))^d$ be the restriction map to the 
subset $\mathfrak{N}_r(X_i)$ (this again is the orthogonal projection if we use
the above convention).

The bounded $U$-equivariant linear map
\begin{equation*}
 \phi\colon l^2(X)^d \to l^2(\mathfrak{N}_r(X_i))^d \oplus \ker(B[i])\colon
 f\mapsto (Q_i(f), \pr_i P_i(f))
\end{equation*}
restricts to a map $\alpha\colon V\to l^2(\mathfrak{N}_r(X_i))^d\oplus
\ker(B[i])$. We claim that $\alpha$ is injective. In fact, if
$f\in V\subset \ker(B_X)$ and $Q_i(f)=0$ then by \eqref{eq:5}
$P_i(f)\in\ker(B[i])$. Therefore, $\pr_i P_i(f)=P_i(f)$. Since
$f\in\ker(\alpha)$ this implies $P_i(f)=0$. But the restriction of
$P_i$ to $V$ is injective, hence $f=0$ as claimed.

It follows that
\begin{equation*}
  \dim_G(\ker(B)) \le
  \frac{1}{\abs{X_i/U}}\dim_U(l^2(\mathfrak{N}_r(X_i))^d\oplus\ker(B[i])).
\end{equation*}
Since $\mathfrak{N}_r(X_i)$ is a free $U$-space, $\dim_U l^2(\mathfrak{N}_r(X_i))^d =
d\abs{\mathfrak{N}_r(X_i)/U}$, therefore
\begin{equation*}
  \dim_G(\ker(B)) \le d\frac{\abs{\mathfrak{N}_r(X_i)/U}}{\abs{X_i/U}} + \frac{\dim_U\ker(B[i])}{\abs{X_i/U}}.
\end{equation*}
By Lemma \ref{lem:amenable_exhaustion}, the first summand on the right 
hand side tends to zero for $i\to\infty$. Consequently
\begin{equation*}
  \dim_G(\ker(B)) \le \liminf_{i\to\infty}
  \frac{\dim_U\ker(B[i])}{\abs{X_i/U}}. 
\end{equation*}
\end{proof}
Using similar ideas, it would be possible to finish the
proof of Theorem \ref{theo:general_amenable_approxi}. However, we will
need more refined estimates on the size of the
spectrum \emph{near} zero to establish Theorem
\ref{newdetest}. Moreover, what follows is in the same way also needed 
for the proof of Theorem \ref{newapproxi}, so we will not finish the
proof now.
Instead, we will continue with the preparation for the proofs of all
the results of Subsections \ref{sec:determinant_bounds} and
\ref{sec:appr-results}, which will comprise the rest of this and the
following three subsections.

We start with two easy observations.

\begin{lemma}\label{lem:trace_normalization}
  With the definitions given in \ref{newsit}, $\tr_i$ in each case is
  a positive and normal trace, which is normalized in the following
  sense: if $\Delta=\id\in M(d\times
  d,\integers G)$ then $\tr_i(\Delta[i])=d$.
\end{lemma}
\begin{proof}
  This follows since $\tr_{G_i}$ has the corresponding
  properties.
\end{proof}

\begin{lemma}\label{ringherit}
Suppose we are in the situation described in \ref{newsit}.
Let $K$ and $L$ be subrings of $\complexs$.
    If $B$ is defined over $KG$,
    then $B[i]$ is defined over $KG_i$ for all $i$ in
    the first and third
    case, and for $i\ge i_0$ in the second case.

   Let
    $\sigma\colon K\to L$ be a ring homomorphism.
Then it induces homomorphisms
    from matrix rings over $KG$ or $KG[i]$ to matrix rings over
$LG$ or $LG[i]$ respectively;
we shall also indicate by $\sigma$ any one of these homomorphisms.
    Furthermore we get
    $\sigma(B[i])=(\sigma B)[i]$ and $B[i]^*=(B^*)[i]$ whenever $B[i]$
    is defined. In particular, the two definitions of $\Delta[i]$ agree.
\end{lemma}
 Since we are
    only interested in $B[i]$ for large $i$, without loss of
    generality, we assume that the statements of Lemma \ref{ringherit}
    are always fulfilled.
\begin{proof}
  The construction of $A[i]$ involve only algebraic (ring)-operations of the
  coefficients of the group elements, and every $\sigma$ is a
  ring-homomorphism, which implies the first and second
    assertion. Taking the adjoint is a
    purely algebraic operation, as well: if $B=(b_{ij})$ then
    $B^*=(b_{ji}^*)$, where for $b=\sum_{g\in G}\lambda_g g\in
    \complexs G$ we have $b^*=\sum_{g\in G}
    \overline{\lambda_{g^{-1}}} g$. Since these algebraic operations
  in the construction of $B[i]$ commute with 
  complex conjugation, the statement about the adjoint operators follows.
\end{proof}

\subsection{An upper bound for the norm of a matrix}
\label{sec:bound}

The approximation results we want to prove here have a sequence of
predecessors in \cite{Lueck(1994c)}, \cite{Dodziuk-Mathai(1998)}, and
\cite{Schick(1998a)}. The first two rely on certain estimates of
operator norms (which are of quite different type in the two cases). 
We want to use a similar approach here, but have to find bounds which
work in both the settings considered in \cite{Lueck(1994c)} and
\cite{Dodziuk-Mathai(1998)}. (These problems were circumvented in
\cite{Schick(1998a)} using a different idea, which however we
didn't manage to apply here). The following definition will supply us
with such a bound.

\begin{definition}\label{def:kappa}
  Let $I$ be an index set. For $A:= (a_{ij})_{i,j\in I}$ with
  $a_{ij}\in\complexs$ set
  \begin{equation*}
S(A):=\sup_{i\in I}\abs{\supp(z_i)},
\end{equation*}
where $z_i$ is the vector
  $z_i:=(a_{ij})_{j\in I}$ and $\supp(z_i):=\{j\in I|\; a_{ij}\ne
  0\}$. Set 
  \begin{equation*}
\abs{A}_{\infty}:=\sup_{i,j}\abs{a_{ij}}.
\end{equation*}
Set $A^*:=(\overline
  a_{ji})_{i,j\in I}$.
If
  $S(A)+S(A^*)+\abs{A}_\infty<\infty$ set 
  \begin{equation*}
\kappa(A):=\sqrt{S(A)S(A^*)}\cdot
  \abs{A}_\infty.
\end{equation*}
Else set $\kappa(A):=\infty$.

  If $G$ is a discrete group and $A\in M(d\times d,\complexs G)$
  consider $A$ as a matrix with complex entries indexed by $I\times I$
  with $I:=\{1,\dots,d\}\times G$, and define $\kappa(A)$ using this
  interpretation.
\end{definition}

\begin{definition}
  Let $I$ and $J$ be two index sets. 
Two matrices $A=(a_{is})_{i,s\in I}$ and $B=(b_{jt})_{j,t\in J}$ are
called \emph{of the same shape}, if for every row or column of
$A$ there exists a
row or column, respectively, of $B$ with the same number of non-zero
entries, and if the sets of non-zero entries of $A$ and $B$ coincide, i.e.
\begin{equation*}
  \{0\}\cup \{ a_{is}\mid i, s\in I\} = \{0\} \cup \{b_{jt}\mid j,t\in J\}.
\end{equation*}
\end{definition}

\begin{lemma}\label{lem:similar_matrices}
  If two matrices $A$ and $B$ have the same shape, then
  \begin{equation*}
    \kappa(A)=\kappa(B).
  \end{equation*}
\end{lemma}
\begin{proof}
  This follows immediately from the definitions.
\end{proof}

\begin{lemma}
  If $A\in M(d\times d,\complexs G)$ as above then $\kappa(A)<\infty$.
\end{lemma}
\begin{proof}
  If $\sum \alpha_gg\in\complexs G$ then $(\sum\alpha_g
  g)^*=\sum\overline{\alpha_{g^{-1}}} g$, and correspondingly for
  matrices.
  
  The assertion follows from the finite support condition and from
  $G$-equivariance $a_{(k,g),(l,h)}= a_{(k,gu),(l,hu)}$ for $g,h,u\in G$.
\end{proof}

\begin{lemma}\label{lem:normest}
  Assume $A=(a_{ij})_{i,j\in I}$ fulfills $\kappa(A)<\infty$. Then $A$
  induces a bounded operator on $l^2(I)$ and for the norm we get
  \begin{equation*}
    \norm{A}\le \kappa(A).
  \end{equation*}
  This applies in particular to $A \in M(d\times d,\complexs G)$
  acting on $l^2(G)^d$.
\end{lemma}
\begin{proof}
  Assume $v=(v_j)\in l^2(I)$. With $z_i:=(a_{ij})_{j\in I}$ we get
  \begin{equation*}
\abs{Av}^2= \sum_{i\in I}\abs{\innerprod{z_i,v}_{l^2(I)}}^2.
\end{equation*}
We estimate
  \begin{equation*}
    \begin{split}
      \abs{\innerprod{z_i,v}}^2 &= \abs{\sum_{j\in\supp(z_i)}
      a_{ij}\overline v_j}^2 \le \sum_{j\in\supp(z_i)}\abs{a_{ij}}^2
    \cdot \sum_{j\in\supp(z_i)}\abs{v_j}^2\\
      &\le S(A) \sup_{i,j\in I}\abs{a_{ij}}^2\cdot \sum_{j\in\supp(z_i)}\abs{v_j}^2.
    \end{split}
  \end{equation*}
  Here we used the Cauchy-Schwarz inequality, but took the size of
the support into account. It follows that
  \begin{equation}\label{Avnorm}
    \abs{Av}^2 \le S(A)\sup_{i,j\in I}\abs{a_{ij}}^2\cdot\sum_{i\in I}\sum_{j\in\supp(z_i)}\abs{v_j}^2.
  \end{equation}
  Observe that
  \begin{equation*}
    j\in \supp(z_i) \iff a_{ij}\ne 0\iff \overline a_{ij}\ne 0 \iff i\in\supp(y_j),
  \end{equation*}
  where $y_j$ is the $j$-th row of $A^*$.
  Consequently, for each fixed $j\in I$, there are not more than
  $S(A^*)$ elements $i\in I$ with $j\in \supp(z_i)$, and hence
  \begin{equation}\label{vnorm}
    \sum_{i\in I}\sum_{j\in \supp(z_i)} \abs{v_j}^2 = \sum_{j\in
    I}\sum_{i\in \supp(y_j)}\abs{v_j}^2\le S(A^*)\abs{v}^2.
\end{equation}
Now \eqref{vnorm} and \eqref{Avnorm} give the desired inequality
\begin{equation*}
  \abs{Av}\le \kappa(A)\cdot \abs{v}. \qed
\end{equation*}
\renewcommand{\qed}{}
\end{proof}

\begin{lemma}\label{cutest}
  Assume $J\subset I$ and $P^*\colon l^2(J)\to l^2(I)$ is the induced
  isometric embedding, with adjoint orthogonal projection $P\colon l^2(I)\to
  l^2(J)$. If $A$ is a matrix indexed by $I$, then $PAP^*$ is a matrix
  indexed by $J$, and 
  \begin{equation*}
    \kappa(PAP^*) \le \kappa(A).
  \end{equation*}
\end{lemma}
\begin{proof}
  This is obvious from the definition, since we only remove entries in
  the original matrix $A$ to get $PAP^*$.
\end{proof}

\begin{lemma}\label{subKequal}
  Suppose $U\subgroup G$. For $A\in M(d\times d,\complexs U)$ let
  $i(A)\in M(d\times d,\complexs G)$ be the image of $A$ under the map
  induced by the inclusion $U\to G$. Then
  \begin{equation*}
    \kappa(A)=\kappa(i(A)).
  \end{equation*}
\end{lemma}
\begin{proof}
  Using the fact that $G$ is a free left $U$-set, and a set of
  representatives of the translates of $U$, we see that
  the matrix for $i(A)$ has copies of $A$ on the ``diagonal'' and
  zeros elsewhere (compare the proof of
  \cite[3.1]{Schick(1998a)}). Therefore the two matrices have the same 
  shape and the statement
  follows from Lemma \ref{lem:similar_matrices}.
\end{proof}

\begin{lemma}\label{lem:X_kappa_equal}
  Let $G$ be a generalized amenable extension of $H$, and $X$ a
  corresponding free $G$-$H$-set as in Definition
  \ref{def:generalized_amenable_extension}. To $A\in M(d\times
  d,\complexs G)$ we associate the operator $A_X\colon l^2(X)^d\to
  l^2(X)^d$. It is now natural to consider $A_X$ as a matrix indexed
  by $I\times I$ with $I:=\{1,\dots,d\}\times X$ (similar to $A$,
  where $X$ is replaced by $G$). With this convention $\kappa(A_X)$ is 
  defined. Then
  \begin{equation*}
    \kappa(A_X) = \kappa(A).
  \end{equation*}
\end{lemma}
\begin{proof}
  The argument is exactly the same as for Lemma \ref{subKequal},
  replacing the free $U$-set $G$ there by the free $G$-set $X$ here.
\end{proof}

We can combine these results to get a uniform estimate for
$\kappa(B[i])$, if we are in the situation described in \ref{newsit}.

\begin{lemma}\label{newdeflim}
  Suppose we are in the situation described in \ref{newsit}. Then 
  \begin{equation*}
    \kappa(B[i])\le \kappa(B)
  \end{equation*}
  for $i$ sufficiently large.
\end{lemma}
\begin{proof}
  In the cases \ref{newsit}(1) and \ref{newsit}(2), we actually prove
  that $\kappa(B[i])=\kappa(B)$ for $i$ sufficiently large. Observe that
  only finitely many elements of $G$ occur as coefficients in $A$. Let
  $Y$ be
  the corresponding subset of $G$. If $G$ is the inverse limit of
  $G_i$ and $G\to G_i$ restricted to $Y$ is injective (which is the
  case for $i$ sufficiently big) the claim follows because the
  matrices describing $B$ and $B[i]$ then have the same shape, and we can
  apply Lemma \ref{lem:similar_matrices}. 

  If $G$ is the direct limit of the groups $G_i$ then, since $A[j_0]$
  has finite support, there is $j_1>j_0$ such that for $i\ge j_1$ the
  map $G_i\to G$ is injective if restricted to the support of
  $A[i]$. We arrive at the conclusion as before.

  In the case \ref{newsit}(3), first observe that by Lemma
  \ref{lem:X_kappa_equal} $\kappa(B)=\kappa(B_X)$. Because of Lemma
  \ref{cutest}, on the other hand, $\kappa(B[i])\le \kappa(B_X)$.
\end{proof}

\subsection{Convergence results for the trace}
\label{sec:trace_convergence}

\begin{lemma}
  \label{lem:trace_conv}
  Suppose we are in the situation described in \ref{newsit}. Assume
  that $p(x)\in\complexs[x]$ is a polynomial. Then
  \begin{equation*}
    \tr_G p(B) =\lim_{i\in I} \tr_i p(B[i])
  \end{equation*}
\end{lemma}
\begin{proof}
  For \ref{newsit}(1) and \ref{newsit}(2), this is proved in
  \cite[Lemma 5.5]{Schick(1998a)}. For \ref{newsit}(3), this is
  essentially the
  content of \cite[Lemma 4.6]{Schick(1998a)}. Actually, there it is
  proved that $\tr_G p(B_X) = \lim_{i\to\infty} \tr_i p(B[i])$, but
  observe that by Subsection \ref{sec:free_actions}, $\tr_G
  p(B_X)=\tr_G p(B)$. Moreover, the assumption that $G$ is a
  generalized amenable extension of $U$ is slightly more general than the
  assumptions made in \cite[Lemma 4.6]{Schick(1998a)}, but
  our assumptions are exactly what is needed in the proof given there.
\end{proof}

\subsection{Kazhdan's inequality}
\label{sec:Kazhdan}

We continue to adopt the situation described in \ref{newsit}.

\begin{definition}
  Define
  \begin{equation*}
    \begin{split}
       \overline{F_\Delta}(\lambda) &:=
       \limsup_{i}F_{\Delta[i]}(\lambda) ,\\
      \underline{F_\Delta}(\lambda) &:=
      \liminf_{i}F_{\Delta[i]}(\lambda) .
    \end{split}
  \end{equation*}
Recall that $\limsup_{i\in I}\{x_i\}=\inf_{i\in I} \{\sup_{j>i}\{x_j\}\}$.
\end{definition}

\begin{definition}
  Suppose $F:[0,\infty)\to \R$ is monotone increasing (e.g.\ a spectral
  density function). Then set
  \begin{equation*}
    F^+(\lambda):= \lim_{\epsilon\to 0^+} F(\lambda+\epsilon)
  \end{equation*}
  i.e.\ $F^+$ is the right continuous approximation of $F$. In
  particular, we have defined $\overline{F_\Delta}^+$ and
  $\underline{F_\Delta}^+$.
\end{definition}

\begin{remark}
  Note that
  by our definition a spectral density function is right continuous,
  i.e.~unchanged if we perform this construction.
\end{remark}

We need the
following functional analytical lemma (compare \cite{Lueck(1994c)} or
\cite{Clair(1999)}):
\begin{lemma}\label{trIE}
  Let $\NeumannN$ be a finite von Neumann algebra with positive normal
  and normalized trace
  $\tr_\NeumannN$. Choose $\Delta\in M(d\times d,\NeumannN)$ positive
  and self-adjoint.\\
  If for a function $p_n:\R\to\R$
\begin{equation}\label{polIE}
 \chi_{[0,\lambda]}(x)\le p_n(x)\le 
  \frac{1}{n}\chi_{[0,K]}(x)+\chi_{[0,\lambda+1/n]}(x) \qquad\forall
  0\le x\le K
\end{equation}
 and if $\norm{\Delta}\le K$ then
 \begin{equation*}
   F_\Delta(\lambda)\le \tr_\NeumannN p_n(\Delta) \le
   \frac{1}{n}d + F_\Delta(\lambda+1/n) .
 \end{equation*}
  Here $\chi_{S}(x)$ is the characteristic function of the subset
  $S\subset \R$.
\end{lemma}
\begin{proof}
  This is a direct consequence of positivity of the trace, of the
  definition of spectral density functions and of the
  fact that $\tr_\NeumannN(1\in M(d\times d,\NeumannN))=d$ by the
  definition of a normalized trace.
\end{proof}

\begin{proposition}\label{appr}
  For every $\lambda\in\reals$ we have
  \begin{equation*} 
    \begin{split}
      &\overline{F_\Delta}(\lambda)\le F_\Delta(\lambda)=
      \underline{F_\Delta}^+(\lambda) ,\\
      &
      F_\Delta(\lambda)
      =\underline{F_\Delta}^+(\lambda)=\overline{F_\Delta}^+(\lambda).
  \end{split}
  \end{equation*}
\end{proposition}
\begin{proof}
  The proof only depends on
     the key lemmata \ref{newdeflim} and \ref{lem:trace_conv}. These
     say (because of Lemma \ref{lem:normest})
    \begin{itemize}
    \item  $\norm{\Delta[i]}\le \kappa(\Delta[i])\le \kappa(\Delta)$ $\forall
      i\in I$
    \item For every polynomial $p\in\complexs[x]$ we have
      $\tr_G(p(\Delta))=\lim_i\tr_i(p(\Delta[i]))$.
    \end{itemize}

  For each $\lambda\in\reals$ choose polynomials $p_n\in\reals[x]$ such that
  inequality \eqref{polIE} is fulfilled. Note that by the first key
  lemma we find the uniform upper bound $\kappa(\Delta)$ for the
  spectrum of all of
  the $\Delta[i]$. Then by Lemma \ref{polIE}
  \begin{equation*}
    \begin{split}
      F_{\Delta[i]}(\lambda) & \le \tr_i( p_n(\Delta[i])) \le
      F_{\Delta[i]}(\lambda+\frac{1}{n}) + \frac{d}{n}
  \end{split}
  \end{equation*}
  We can take $\liminf$ and $\limsup$
  and use the second key lemma to get
  \begin{equation*}
    \overline{F_\Delta}(\lambda) \le \tr_G(p_n(\Delta)) \le
    \underline{F_\Delta}(\lambda+\frac{1}{n}) +\frac{d}{n}.
  \end{equation*}
  Now we take the limit as $n\to\infty$. We use the fact that
  $\tr_G$ is normal and $p_n(\Delta)$ converges strongly inside a
  norm bounded set to $\chi_{[0,\lambda]}(\Delta)$. Therefore the
  convergence even is in the ultra-strong topology.

  Thus
  \begin{equation*}
    \overline{F_\Delta}(\lambda) \le F_\Delta(\lambda) \le
    \underline{F_\Delta}^+(\lambda) .
  \end{equation*}
  For $\epsilon>0$ we can now conclude, since $\underline{F_\Delta}$
  and $\overline{F_\Delta}$ are monotone, that
\begin{equation*}
 F_\Delta(\lambda)\le \underline{F_\Delta}(\lambda+\epsilon)\le
 \overline{F_\Delta}(\lambda+\epsilon)\le F_\Delta(\lambda+\epsilon).
\end{equation*}
Taking the limit as $\epsilon\to 0^+$ gives (since $F_\Delta$ is right
continuous)
\begin{equation*}
F_\Delta(\lambda)=\overline{F_\Delta}^+(\lambda)=
\underline{F_\Delta}^+(\lambda).
\end{equation*}
Therefore both inequalities are established.
\end{proof}

We are now able to finish the proof of Theorem
\ref{theo:general_amenable_approxi}. 
\begin{proof}[Proof of Theorem \ref{theo:general_amenable_approxi}]
  We specialize Proposition \ref{appr} to $\lambda=0$  and see that
  \begin{equation*}
    \dim_G(\ker(\Delta))\ge \limsup_{i} \dim_i(\ker(\Delta[i])).
  \end{equation*}
  On the other hand, by Lemma \ref{lem:half_proof},
  \begin{equation*}
      \dim_G(\ker(\Delta)) \le \liminf_{i}
      \dim_i\ker(\Delta[i]) ,
\end{equation*}
and consequently
\begin{equation*}
      \dim_G(\ker(\Delta)) = \lim_{i}
  \dim_i\ker(\Delta[i]).
\end{equation*}
  As remarked above, the conclusion for an arbitrary $B$ follows by
  considering $\Delta=B^*B$.
\end{proof}

To prove an estimate similar to Lemma \ref{lem:half_proof} in the
general situation of \ref{newsit}, we need additional control on the
spectral measure near zero. This can be given using uniform bounds 
on determinants.

\begin{theorem}\label{theo:strong_convergence}
  In the situation described in \ref{newsit}, assume that there is
  $C\in\reals$ such that
  \begin{equation*}
    \ln\det_i(\Delta[i]) \ge C\qquad\forall i\in I.
  \end{equation*}
  Then $\ln\det_G(\Delta)\ge C$, and
  \begin{equation*}
    \dim_G(\ker(\Delta)) =\lim_{i\in I} \dim_i(\ker(\Delta[i])).
  \end{equation*}
\end{theorem}
\begin{proof}
   Set $K:=\kappa(\Delta)$. Then, by Lemma \ref{lem:normest} and Lemma 
   \ref{newdeflim}, $K>\norm{\Delta}$ and $K>\norm{\Delta[i]}$ $\forall
 i$. Hence,
  \begin{equation*}
    \ln\Det_i(\Delta[i])=\ln(K)(F_{\Delta[i]}(K)-F_{\Delta[i]}(0))
    -\int_{0^+}^{K}
    \frac{F_{\Delta[i]}(\lambda)-F_{\Delta[i]}(0)}{\lambda}\; d\lambda .
  \end{equation*}
  If this is (by assumption) $\ge C$, then since
  $F_{\Delta[i]}(K)=\tr_i(\id_d)=d$ by our normalization
  \begin{equation*}
    \int_{0^+}^{K}
    \frac{F_{\Delta[i]}(\lambda)-F_{\Delta[i]}(0)}{\lambda} \;d\lambda
    +C \le
    \ln(K)(d-F_{\Delta[i]}(0))\le \ln(K)d .
  \end{equation*}
 We want to establish the same estimate for
  $\Delta$. If $\epsilon>0$ then
  \begin{equation*}
    \begin{split}
       \int_{\epsilon}^K\frac{F_\Delta(\lambda)-F_\Delta(0)}{\lambda}\;d\lambda
      & = \int_{\epsilon}^K
      \frac{\underline{F_\Delta}^+(\lambda)-F_\Delta(0)}{\lambda}\;d\lambda
      =\int_{\epsilon}^K\frac{\underline{F_\Delta}(\lambda)-F_\Delta(0)}{\lambda}\\
      \intertext{(since  the integrand is bounded, the integral
        over the left continuous approximation is equal to the integral
        over the original function)}\\ 
      \le & \int_{\epsilon}^K
      \frac{\underline{F_\Delta}(\lambda)-\overline{F_\Delta}(0)}{\lambda}\\
      = & \int_\epsilon^K \frac{\liminf_i F_{\Delta[i]}(\lambda)-\limsup_i
        F_{\Delta[i]}(0)}{\lambda}\\
      \le &\int_\epsilon^K
      \frac{\liminf_i (F_{\Delta[i]}(\lambda)-F_{\Delta[i]}(0))}{\lambda}\\
      \le &\liminf_i\int_\epsilon^K
      \frac{F_{\Delta[i]}(\lambda)-F_{\Delta[i]}(0)}{\lambda} \le
      d\ln(K)-C .
  \end{split}
\end{equation*}
Since this holds for every $\epsilon>0$, we even have
\begin{equation*}
  \begin{split}
    \int_{0^+}^K \frac{F_\Delta(\lambda)-F_\Delta(0)}{\lambda} 
    \le & \int_{0^+}^K \frac{\underline{F_\Delta}(\lambda)-
      \overline{F_\Delta}(0)}{\lambda}\;d\lambda\\
    \le & \sup_{\epsilon>0}\liminf_{i} \int_{\epsilon}^K
    \frac{F_{\Delta[i]}(\lambda)-F_{\Delta[i]}(0)}{\lambda}\;d\lambda\\
    \le & d \ln(K) -C.
  \end{split}
\end{equation*}
The second integral would be infinite if $\lim_{\delta\to
  0}\underline{F_\Delta}(\delta)\ne \overline{F_\Delta(0)}$. It
follows that $\limsup_i
F_{\Delta[i]}(0)=F_\Delta(0)$. Since we can play the same game for every
subnet of $I$, also $\liminf_i F_{\Delta[i]}(0)=F_\Delta(0)$ i.e.\ the
approximation result is true.

For the estimate of the determinant note that in the above inequality 
\begin{equation*}
  \begin{split}
    & \sup_{\epsilon>0}\liminf_{i}\int_{\epsilon}^K
    \frac{F_{\Delta[i]}(\lambda)-F_{\Delta[i]}(0)}{\lambda}\;d\lambda\\
    \le & \liminf_i \sup_{\epsilon>0}\int_\epsilon^K
    \frac{F_{\Delta[i]}(\lambda)-F_{\Delta[i]}(0)}{\lambda}\;d\lambda = \liminf_i \int_{0^+}^K
    \frac{F_{\Delta[i]}(\lambda)-F_{\Delta[i]}(0)}{\lambda}\;d\lambda\\
    \le & \ln(K)(d-F_{\Delta[i]}(0))-C .
  \end{split}
\end{equation*}
Therefore
\begin{equation*}
  \begin{split}
    \ln\Det_G(\Delta) &=\ln(K)(d-F_\Delta(0)) -\int_{0^+}^K
    \frac{F_\Delta(\lambda)-F_\Delta(0)}{\lambda}\;d\lambda\\
    &\ge \ln(K)(d-\lim_i F_{\Delta[i]}(0)) - \liminf_i \int_{0^+}^K
    \frac{F_{\Delta[i]}(\lambda)-F_{\Delta[i]}(0)}{\lambda}\;d\lambda\\
    &= \limsup_i\left(\ln(K)(d-F_{\Delta[i]}(0)) - \int_{0^+}^K
      \frac{F_{\Delta[i]}(\lambda)-F_{\Delta[i]}(0)}{\lambda}\;d\lambda\right)\\
    & 
\ge C. \qquad\qed
  \end{split}
\end{equation*}
\renewcommand{\qed}{}
\end{proof}

\subsection{Proof of Theorems \protect\ref{newdetest} and
\protect\ref{newapproxi} by induction}

Theorems \ref{newdetest} and \ref{newapproxi} are generalizations of
\cite[Theorem 6.9]{Schick(1998a)}, where the corresponding statements
are proved for matrices over $\integers G$. The proof will be done by
transfinite induction (using the induction principle
\cite[2.2]{Schick(1998a)}) and is very similar to the proof of
\cite[6.9]{Schick(1998a)}. Therefore, we will not give all the details
here but concentrate on the necessary modifications.

We will use the following induction principle.
  \begin{proposition}\label{ind_princ}
    Suppose a property $C$ of groups is shared by the trivial group, and
    the following is true:
    \begin{itemize}
    \item  whenever $U$ has property $C$ and $G$ is a generalized
      amenable extension of $U$, then $G$ has property $C$ as well;
    \item whenever $G$ is a direct or inverse limit of a directed
      system of groups with property $C$, then $G$ has
      property $C$
    \item If $U\subgroup G$, and $U$ has property $C$, then also $G$
      has property $C$.
    \end{itemize}
    Then property $C$ is shared by all groups in
    the class $\RAgroups$.
  \end{proposition}
  \begin{proof}
    The proof of the induction principle is done by transfinite
    induction in a standard way, using the definition of
    $\RAgroups$. The result corresponds to \cite[Proposition
    2.2]{Schick(1998a)}. 
\end{proof}

We are going to use the induction principle to prove the determinant
bound property and the algebraic continuity property of Definition
\ref{def:detestprop}.

\subsubsection{Trivial group}

First we explain how the induction gets started:
\begin{lemma}
  The trivial group $G=\{1\}$ has the determinant bound property and
  the algebraic continuity property.
\end{lemma}
\begin{proof}
  Remember that $\Delta=B^*B$ and $\det_{\{1\}}(B^*B)$ in this case is the
  product of
all non-zero eigenvalues of $B^*B$, which is positive since $B^*B$ is
a non-negative operator. If $q(t)$ is the characteristic
polynomial of $\Delta$ and $q(t)=t^\nu \overline q(t)$ with $\overline
q(0)\ne 0$ then $\det_{\{1\}}(\Delta)=\overline q(0)$. In particular it
is contained in $o(\overline\rationals)$. Moreover
$\det_{\{1\}}(\sigma_k(A))=\sigma_k(\det_{\{1\}}(A))$ for $ k=1,\dots,r$. The
product of all these numbers is fixed under $\sigma_1,\dots,\sigma_r$,
therefore a rational number and is an algebraic integer, hence a (non-zero)
integer. Taking absolute values, we get a positive integer, and consequently
\begin{equation*}
  \sum_{k=1}^r \ln\abs{\det_{\{1\}}({\sigma_k(A)})} \ge 0.
\end{equation*}
Using Lemma \ref{q0approxi} we get the desired estimate of the
determinant.

Since $\complexs$ is flat over $L$,
\begin{equation*}
\dim_{\{1\}}(\ker(B))=\dim_\complexs(\ker(B))=\dim_L(\ker(B|_{L^d})).
\end{equation*}
Choose $\sigma\in\{\sigma_1,\dots,\sigma_r\}$. Since
$\sigma\colon L\to\sigma(L)$ is a field automorphism,
\begin{equation*}
\dim_L(\ker(B|_{L^d}))=\dim_{\sigma(L)}(\ker(\sigma(B)|_{\sigma(L)^d})),
\end{equation*}
which by the same
reasoning as above coincides with $\dim_{\{1\}}(\ker(\sigma B))$. 
\end{proof}

\begin{lemma}\label{q0approxi}
  Let $A\in M(d\times d,\complexs)$. Then
  \begin{equation*}
    \abs{ \det_{\{1\}}(A) } \le \norm{A}^d,
  \end{equation*}
  where $\norm{A}$ denotes the Euclidean operator norm.
\end{lemma}
\begin{proof}
  The result is well known if $A$ has no kernel and
  $\det_{\{1\}}(A)=\det(A)$. By scaling, we may assume that
  $\norm{A}\ge 1$. If $A$ has a non-trivial kernel, choose a basis
  whose first elements span the kernel of $A$. Then $A=
  \begin{pmatrix}
    0 & * \\ 0 & A_0
  \end{pmatrix}$ and $A_0$ has a smaller dimension than $A$. It is
  obvious from the definitions that
  $q_A(t)=t^{\dim(\ker(A))}q_{A_0}(t)$, where $q_A$ is the
  characteristic polynomial of $A$. Therefore
  $\det_{\{1\}}(A)=\det_{\{1\}}(A_0)$ (since this is the first
  non-trivial coefficient in the characteristic polynomial). By
  induction $\det_{\{1\}}(A)\le \norm{A_0}^{d-\dim{\ker(A)}}$. Since
  $\norm{A_0}\le \norm{A}$, the result follows. Observe that the
  statement is trivial if $d=1$.
\end{proof}

\subsubsection{Subgroups}

We have to
check that, if $G$ has the determinant bound property and the
algebraic continuity
property, then the same is true for any subgroup $U$ of $G$.

 This is done as
follows: a matrix over $\complexs U$ can be induced up to (in other
words, viewed
as) a matrix
over $\complexs G$. By Lemma \ref{subKequal}, the right hand side of
\eqref{basicest} is
unchanged under this process, and it is well known that the spectrum
of the operator and therefore \eqref{sigmaequal} and the left hand
side of \eqref{basicest} is unchanged, as
well, compare e.g.~\cite[3.1]{Schick(1998a)}. Since the spectrum is
unchanged, the algebraic continuity also is immediately inherited.

\subsubsection{Limits and amenable extensions}

Now assume $G$ is a generalized amenable extension of $U$, or the direct or
inverse limit of a directed system of groups $(G_i)_{i\in I}$, and
assume $B\in M(d\times d,o(\overline\rationals) G)$, and
$\Delta=B^*B$. That means, we are exactly in the situation described
in \ref{newsit} (as we have seen before, without loss of generality we 
can assume throughout that the coefficients are algebraic integers
instead of
algebraic numbers).

We want to use the induction hypothesis to establish a uniform lower
bound on $\ln\det_i(\Delta[i])$.
\begin{lemma}
  In the given situation,
  \begin{equation}\label{eq:uniform_det_est}
    \ln\det_i(\Delta[i]) \ge  -d \sum_{k=2}^r
    \ln(\kappa(\sigma_k(\Delta))) \quad\text{for $i$ sufficiently large.}
  \end{equation}
\end{lemma}
\begin{proof}
  In case of \ref{newsit}(1) or \ref{newsit}(2), since
  $\det_i=\det_{G_i}$ in this situation, by the induction hypothesis
  we have
  \begin{equation*}
    \ln\det_i(\Delta[i]) \ge  -d \sum_{k=2}^r
    \ln(\kappa(\sigma_k(\Delta[i]))).
  \end{equation*}
  However, by Lemma \ref{ringherit},
  $\sigma_k(\Delta[i])=\sigma_k(\Delta)[i]$ if $i$ is sufficiently
  large, moreover, $\kappa(\sigma_k(\Delta)[i])\le
  \kappa(\sigma_k(\Delta))$ by Lemma \ref{newdeflim}, again for $i$
  sufficiently large. The inequality \eqref{eq:uniform_det_est} follows.

  It remains to treat the case \ref{newsit}(3). Here,
  \begin{equation}
\ln\det_i(\Delta[i]) =\frac{1}{N_i}\ln\det_U(\Delta[i]).
\end{equation}
On the other hand, $\Delta[i]$ is to be considered not as an element
of $M(d\times d, \complexs U)$, but as an element of $M(N_id\times N_i 
d,\complexs U)$. Consequently, by the induction hypothesis,
\begin{equation*}
  \ln\det_U(\Delta[i]) \ge -N_i d\sum_{k=2}^r
    \ln(\kappa(\sigma_k(\Delta[i]))).
  \end{equation*}
  Dividing this inequality by $N_i>0$, and using the fact that by
  Lemma \ref{newdeflim}
  $\kappa(\sigma_k(\Delta))\ge \kappa(\sigma_k(\Delta[i]))$, again we
  conclude that the inequality \eqref{eq:uniform_det_est} is true.
\end{proof}

Because of Theorem \ref{theo:strong_convergence},
\begin{equation}
    \ln\det_G(\Delta) \ge  -d \sum_{k=2}^r
    \ln(\kappa(\sigma_k(\Delta))).
\end{equation}
The induction principle of Proposition \ref{ind_princ} therefore implies
that, if
$G\in\RAgroups$, then $G$ has the determinant bound property.

Note, moreover, that we can now prove Theorem \ref{newapproxi} inductively.

To carry out the induction step, we note that
\begin{equation*}
  \ker(\sigma_kB)=\ker((\sigma_kB)^*\sigma_k B)\qquad k=1,\dots,r.
\end{equation*}
The approximation result we have just proved implies (using Lemma
\ref{ringherit}) that 
\begin{equation*}
  \dim_G(\ker((\sigma_k B)^*\sigma_k B)) = \lim_{i\in I}
  \dim_i(\ker((\sigma_k B)[i]^*\sigma_k B[i])).
\end{equation*}
Again we know that for each $i$ $\ker((\sigma_k
B)[i]^*\sigma_kB[i])=\ker(\sigma_k B[i])$, and the induction
hypothesis implies that $\dim_i\ker(\sigma_k B[i])=\dim_i\ker(B[i])$ for
each $k=1,\dots,r$ and for each $i$. Therefore, we obtain
\begin{equation*}
  \dim_G(\ker(B)) = \dim_G(\ker(\sigma_k(B))).
\end{equation*}
By induction, if $G\in\RAgroups$ then $G$ has the algebraic
continuity property.

This finishes the proof of Theorems \ref{newdetest} and \ref{newapproxi}.

\subsection{Proof of Proposition \protect\ref{prop:algebraiclimits}}

To conclude this section, we note that
Theorem \ref{newapproxi} immediately implies Proposition
\ref{prop:algebraiclimits} because the $L^2$-dimension of the kernel we
have to compute is the limit of $L^2$-dimensions which by assumption
are integers. Since $\integers$ is discrete, the limit has to be an
integer, too.

\section{Absence of transcendental eigenvalues}

For a discrete amenable group $G$, consider self-adjoint operators
$A=B^\ast B$, $B\in M(d\times d,\complexs G)$,
regarded as operators acting on the left on $l^2(G)^d$.
Then as a consequence of Theorem \ref{theo:general_amenable_approxi}
and Proposition \ref{appr},
\begin{gather}
  \label{eq:amenablespecincl}
  \spec A \subset \overline{\bigcup_{i=1}^\infty \spec A[i]}\\
  \label{eq:amenableptspec}
  \specpt A \subset \bigcup_{i=1}^\infty \spec A[i],
\end{gather}
where the $A[i]$ (as in \ref{newsit} part 3 with group $U$ trivial)
are finite dimensional matrices over $\complexs$.
These results have also been shown for the particular case when $A$ is
the Laplacian in \cite{Mathai-Yates(2002)}.

This raises the question: do analogues of
(\ref{eq:amenablespecincl}) and (\ref{eq:amenableptspec}) hold for
any non-amenable groups? When the
group algebra is restricted to $\overline\rationals G$, the
answer is yes: we show by Corollary \ref{cor:aepspec}
that (\ref{eq:amenablespecincl})
holds for any group which has the algebraic eigenvalue property,
defined below. Such groups include the amenable groups and additionally
any group in the class $\LinnelsC$, which includes
the free groups and is closed under extensions with elementary amenable quotient
and under directed unions.
If further the groups $G_i$ (as in \ref{newsit}) all belong to the
class $\RAgroups$, Corollary \ref{cor:aepdimeq} demonstrates that
the point spectrum inclusion (\ref{eq:amenableptspec}) also holds.
A weaker result holds for
groups $G$ in the larger class $\RAgroups$, namely that a self-adjoint
$A\in M(d\times d,\overline\rationals G)$ has no eigenvalues that are Liouville
transcendental numbers (Theorem \ref{theo:non_eigenvalue}).

A consequence of (\ref{eq:amenableptspec}) is that for
an amenable group $G$ and $A=B^\ast B$,
$B\in M(d\times d,\overline\rationals G)$, the point spectrum
of $A$ is a subset of the algebraic numbers. This prompts the
definition of the algebraic eigenvalue property as follows.

\begin{definition}
  \label{algeigenprop}
  We say that a discrete group $G$ has the \emph{algebraic eigenvalue
    property} (or the \emph{algebraic eigenvalue property for
    algebraic matrices}), if for every matrix $A\in M(d\times d,
  \overline\rationals G)$ the eigenvalues of $A$, acting on $l^2(G)^d$
  are algebraic. We say $G$ has the \emph{algebraic eigenvalue
    property for rational matrices}, if the same is true for every
  matrix $A\in M(d\times d,\rationals G)$.
\end{definition}

Observe that for a group $G$ with the algebraic eigenvalue property, and
a matrix $A$ over $\overline\rationals G$, a
transcendental number is necessarily a point of continuity for the
spectral density function of $A$. 

\begin{corollary}
  \label{cor:aepdimeq}
  Assume that we are in the situation described in \ref{newsit}, and 
  that $G$
  has the algebraic eigenvalue property. If $B\in M(d\times d,
  \overline\rationals G)$ and $\lambda\in\mathbb{C}$ is a
  transcendental number, then 
  \begin{equation}\label{eq:3}
    \dim_G(\ker(B-\lambda)) = \lim_{i\in I}\dim_i(\ker(B[i]-\lambda)).
  \end{equation}
  Moreover, if $G_i\in \RAgroups$, then \eqref{eq:3} holds for every
  number $\lambda\in\mathbb{C}$.
\end{corollary}
\begin{proof}
We shall repeatedly use the fact that $\ker (A^*A) = \ker A$ for a
general matrix $A$.
Set $\Delta = (B-\lambda)^*(B-\lambda)$.  Then we have (cf.\ Lemma
\ref{ringherit}) $\Delta[i] =
(B[i]-\lambda)^*(B[i]-\lambda)$ (for $i \ge i_0$).  If $\lambda$ is
transcendental, then by the very definition of the algebraic
eigenvalue property, we see that $\lambda$ is not an eigenvalue of
$B$.  Therefore $\ker(B-\lambda) = 0$ and
consequently $\ker \Delta = 0$.  We now deduce from
\cite[Lemma 7.1]{Schick(1998a)} that
$\lim_{i\in I}\dim_i(\ker(\Delta[i])) = 0$ and we conclude that
$\lim_{i\in I}\dim_i(\ker(B[i] - \lambda)) = 0$.  This proves the
first statement and the second statement in the case $\lambda$ is
transcendental.  Finally when $\lambda$ is algebraic,
the second statement follows from Theorem \ref{newapproxi}.
\end{proof}


\begin{corollary}
  \label{cor:aepspec}
  Let $B\in M(d\times d,\overline\rationals G)$, where $G$ has
  the algebraic eigenvalue property, as in Corollary \ref{cor:aepdimeq}.
  Following the situation \ref{newsit}, take $A=B^\ast B$
  and $A[i]=B[i]^\ast B[i]$.
  Then regarding $A$ as an operator on $l^2(G)^d$,
  \begin{equation}
    \label{eq:specincl}
    \spec A \subset
    \overline{\bigcup_{i\in I} \spec A[i]}.
  \end{equation}
\end{corollary}
\begin{proof}
  Suppose $\lambda_1$,$\lambda_2$ are transcendental real numbers,
  with
  \[
  [\lambda_1,\lambda_2] \cap \bigcup_{i\in I} \spec A[i]
  =\emptyset .
  \]
  By the definition of the algebraic eigenvalue property,
  the spectral density function $F_A$ is continuous at
  transcendental $\lambda$, and so by Proposition \ref{appr}
  \[
  F_A(\lambda)=\lim_{i\in I} F_{A[i]}(\lambda)
  \qquad\text{for $\lambda$ transcendental in $\reals$.}
  \]
  As $[\lambda_1,\lambda_2]$ is in a gap of the spectrum of every
  $A[i]$,
  \[
    F_A(\lambda_2)-F_A(\lambda_1)
    = \lim_{i\in I} F_{A[i]}(\lambda_2)-F_{A[i]}(\lambda_1)
    =0.
  \]
  Therefore
  \begin{equation}
    \label{eq:gapinclusion}
    [\lambda_1,\lambda_2] \cap \bigcup_{i\in I} \spec A[i]
    =\emptyset
    \implies
    (\lambda_1,\lambda_2) \cap \overline{\spec A}
    =\emptyset .
  \end{equation}

  Consider a point $p$ in $\spec A$, and sequences
  $\lambda^+_j$ and $\lambda^-_j$ of real transcendental numbers
  that converge to $p$ from above and below respectively. Then by
  (\ref{eq:gapinclusion}),
  \[
  [\lambda^-_j,\lambda^+_j] \cap
  \overline{\bigcup_{i\in I} \spec A[i]}
  \neq\emptyset.
  \]
  These constitute a strictly decreasing sequence of closed, non-empty subsets
  of the real line, and so have non-empty intersection.
  \begin{align*}
    \emptyset\neq
    \bigcap_{j=1}^\infty \Big(
    [\lambda^-_j,\lambda^+_j] \cap
    \overline{\bigcup_{i\in I} \spec A[i]}
    \Big)
    &=
    \bigcap_{j=1}^\infty
    [\lambda^-_j,\lambda^+_j]
    \cap \overline{\bigcup_{i\in I} \spec A[i]}
    \\
    &=
    \{ p \} \cap \overline{\bigcup_{i\in I} \spec A[i]}.
  \end{align*}
  As this holds for any $p$ in the spectrum of $A$,
  (\ref{eq:specincl}) follows.
\end{proof}


\begin{example}
  The trivial group has the algebraic eigenvalue property, since the
  eigenvalues are the zeros of the characteristic polynomial. The same 
  is true for every finite group. More generally, if $G$ contains a
  subgroup $H$ of finite index, and $H$ has the algebraic eigenvalue
  property, then the same is true for $G$. And if $G$ has the
  algebraic eigenvalue property and $H$ is a subgroup of $G$, then $H$ 
  has the algebraic eigenvalue property, too.
\end{example}
\begin{proof}
  This follows immediately from \cite[Proposition 3.1]{Schick(1998a)}
  and from \cite[Lemma 8.6]{Linnell(1998)}.
\end{proof}

The following theorem
has already been established by Roman Sauer in the
special case that the relevant matrix $A$ is self-adjoint; see
\cite{Sauer(2001)} for this and many stronger results.

\begin{theorem}
  \label{thm:freegroups_aep}
  Free groups have the algebraic eigenvalue property.
\end{theorem}

We will deduce this from the following result.
\begin{theorem} \label{thm:ordergroups_aep}
Let $G$ be an ordered group and suppose $G$ satisfies the strong
Atiyah conjecture over $\mathbb {C}G$.  Then $G$ has the algebraic
eigenvalue property.
\end{theorem}
To say that $G$ is an ordered group means that it has a total order
$\le$ with the property that if $x\le y$ and $g \in G$, then $xg \le
yg$ and $gx \le gy$.  Since free groups are orderable and $G$
satisfies the strong Atiyah conjecture over $\mathbb {C}G$
\cite[Theorem 1.3]{Linnell(1993)}, Theorem \ref{thm:freegroups_aep}
follows from Theorem \ref{thm:ordergroups_aep}.

\begin{proof}[Proof of Theorem \ref{thm:ordergroups_aep}]
For an arbitrary group $G$, let $\universalU G$ indicate the ring of
operators affiliated to the group von Neumann algebra $\NeumannN G$
\cite[\S 8]{Linnell(1998)} and
\cite[Definition 7, p.~741]{Schick(1999)}; this is a ring containing
$\NeumannN G$ in which every element is either a zero divisor or
invertible.  Also let $D(G)$ denote the division closure of
$\mathbb{C}G$ in $\universalU G$, that is the smallest subring of
$\universalU G$ which contains $\mathbb {C}G$ and is closed under
taking inverses.  Of course if $H \le G$, then $\universalU H$ is
naturally a subring of $\universalU G$ and hence $D(H)$ is also a
subring of $D(G)$.  Furthermore if $x_1, x_2, \dots$ are in distinct
right cosets of $H$ in $G$, then the sum $\sum_i (\universalU
H)x_i$ is direct and in particular the sum $\sum_i D(H)x_i$ is also
direct.

We describe the concept of a free division ring of fractions for
$KG$ where $K$ is a field, as defined on
\cite[p.~182]{Hughes(1970)}.
This is a division ring $D$ containing $KG$ which is generated by
$KG$.  Also if $D_N$ denotes the division closure of
$KN$ in $D$ for the subgroup $N$ of $G$, then it has the
following property.  If $N \lhd H \le G$, $H/N$ is infinite cyclic
and generated by $Nt$ where $t \in H$, then the sum $\sum_i
D_N t^i$ is direct.

Now let $G$ be an ordered group which satisfies the Atiyah conjecture
over $\mathbb {C}G$.  Then $D(G)$ is a division ring by
\cite[Lemma 3]{Schick(1999)}, and by the first paragraph is a
free division ring of fractions for $\mathbb {C}G$.  Also
we can form the
Malcev-Neumann division ring $\mathbb {C}[[G]]$ as described in 
\cite[Corollary 8.7.6]{Cohn(1985)}.  The elements of $\mathbb
{C}[[G]]$ are power series of the form $\sum_{g\in G} a_gg$ with
$a_g \in \mathbb {C}$ whose support $\{g\in G \mid a_g \ne 0\}$ is
well ordered.  If $M$ indicates the division
closure of $\mathbb {C}G$ in $\mathbb {C}[[G]]$, then $M$ is also
a free division ring of fractions for $\mathbb {C}G$,
so by \cite[Theorem,
p.~182]{Hughes(1970)} there is an isomorphism from $D(G)$ onto $M$
which extends the identity map on $\mathbb {C}G$.  Therefore we may
consider $D(G)$ as a subring of $\mathbb {C}[[G]]$.

Let $d$ be a positive integer and let $A \in M(d \times d,
\overline{\mathbb {Q}}G)$.  We shall repeatedly use the following
fact without comment: if $D$ is a division ring and $x \in M(d
\times d, D)$, then $x$ is a zero divisor if and only if it is not
invertible.  Also $x$ is a right or left zero divisor if and only if
it is a two-sided zero divisor, and is left or right invertible if
and only if it is two-sided invertible.
By \cite[Lemma 4.1]{Linnell(1993)}, the division closure of
$M(d\times d, \mathbb {C}G)$ in $M(d\times d, \universalU G)$
is $M(d \times d, D(G))$.  Suppose $A$ has a transcendental
eigenvalue $1/t$.  Then $I - tA$ is not invertible in $M(d\times d,
\universalU G)$, so $I - tA$ is not invertible in $M(d\times d,D(G))$.
Therefore $I - tA$ is a zero divisor in $M(d\times d, D(G))$ and we
deduce that $I - tA$ is not invertible in $M(d\times d, \mathbb
{C}[[G]])$.  Thus $I - tA$ is not invertible in $M(d\times d,
\overline{\mathbb {Q}}(t)[[G]])$ and we conclude that $I - tA$ is a
zero divisor in $M(d\times d, \overline{\mathbb {Q}}(t)[[G]])$.
Let us write $T$ for the infinite cyclic group generated by $t$.
Then $\overline{\mathbb {Q}}(t)$ embeds in
$\overline{\mathbb {Q}}[[T]]$ and we deduce
that $I - tA$ is a zero divisor in $M(d\times d, (\overline{\mathbb
{Q}}[[T]])[[G]])$. Let $E$ denote the division closure of
$\overline{\mathbb {Q}}[T\times G]$ in
$(\overline{\mathbb {Q}}[[T]])[[G]]$.
Then $I-tA$ is not invertible in $M(d\times d,E)$ and hence is a
zero divisor in $M(d\times d, E)$.
Now $(\overline{\mathbb {Q}}[[T]])[[G]]$ is just the
Malcev-Neumann power series ring with
respect to the ordered group $T\times G$, where we have given $T
\times G$ the lexicographic ordering; specifically $(t^i,h) <
(t^j,g)$ means $h<g$ or $h=g$ and $i < j$.  Note that $E$ is a free
division ring of fractions for $\overline{\mathbb {Q}}[T\times G]$.
We may also form the
Malcev-Neumann power series with respect to the ordering $(t^i,h) <
(t^j,g)$ means $i<j$ or $i=j$ and $h<g$.  The power series ring we
obtain in this case is $(\overline{\mathbb {Q}}[[G]])[[T]]$.  Since
the division closure of $\mathbb {Q}[T\times G]$ in
$(\overline{\mathbb {Q}}[[G]])[[T]]$ is also a free division ring of
fractions for $\mathbb {Q}[T\times G]$, we may by \cite[Theorem,
p.~182]{Hughes(1970)} view $E$ as a subring of $(\overline{\mathbb
{Q}}[[G]])[[T]]$.  We deduce that $I-At$ is a zero divisor in
$(\overline{\mathbb {Q}}[[G]])[[T]]$, which we see by a leading
term argument is not the case.
\end{proof}

\begin{theorem}\label{theo:algebraic_eigenvalues_and_amenable_extensions}
  Assume $H$ has the algebraic eigenvalue property, or the algebraic
  eigenvalue property for rational matrices. Let $G $ be a
  generalized amenable extension of $H$. Then $G $ also has the
  algebraic eigenvalue property, or the algebraic eigenvalue property 
  for rational matrices, respectively.
\end{theorem}
\begin{proof}
  Let $B\in M(d\times d,\overline\rationals G)$. By
  Theorem \ref{theo:general_amenable_approxi}, as
  $B-\lambda\in M(d\times d,\complexs G)$,
  \[
  \dim_G\ker (B-\lambda)=\lim_{i\to\infty}
  \frac{\dim_H (B[i]-\lambda)}{N_i}.
  \]
  The $B[i]$ are matrices over $\overline\rationals H$, and so
  if $H$ has the algebraic eigenvalue property,
  $\dim_H B[i]-\lambda=0$ for all transcendental $\lambda$.
  So
  \[
  \lambda\not\in\overline\rationals \implies \dim_G\ker (B-\lambda)=0,
  \]
  that is, $B$ has only algebraic eigenvalues. $G$ therefore
  has the algebraic eigenvalue property. The same argument applies
  for the case where $H$ has the algebraic eigenvalue property
  for rational matrices.
\end{proof}

\begin{corollary}
  \label{Chastheaep}
  Every amenable group, and every group in the class $\LinnelsC$ of
  Linnell has the algebraic eigenvalue property.
\end{corollary}
\begin{proof}
  Since the trivial group has the algebraic eigenvalue property,
  Theorem \ref{theo:algebraic_eigenvalues_and_amenable_extensions}
  implies the statement for amenable groups.

  The proof of the algebraic eigenvalue property for groups in
  $\LinnelsC$ is done by transfinite induction, following the
  pattern of \cite{Linnell(1993)}.

  For ordinals $\alpha$ define the class of groups $\LinnelsC_\alpha$
  as follows:
  \begin{enumerate}
  \item $\LinnelsC_0$ is the class of free groups.
  \item $\LinnelsC_{\alpha+1}$ is the class of groups $G$ such that
    $G$ is the directed union of groups $G_i\in\LinnelsC_\alpha$,
    or $G$ is the extension of a group $H\in\LinnelsC_\alpha$ with
    elementary amenable quotient.
  \item $\LinnelsC_\beta=\bigcup_{\alpha<\beta}\LinnelsC_\alpha$ when
    $\beta$ is a limit ordinal.
  \end{enumerate}
  Then a group is in $\LinnelsC$ if it belongs to $\LinnelsC_\alpha$ for
  some ordinal $\alpha$. The algebraic eigenvalue property holds for
  groups in $\LinnelsC_0$ by Theorem \ref{thm:freegroups_aep}.
  Proceeding by transfinite induction, we need to establish that
  groups in $\LinnelsC_\beta$ have the algebraic eigenvalue property
  given that all groups in $\LinnelsC_\alpha$ have the property
  for all $\alpha<\beta$. For limit ordinals $\beta$, this follows
  trivially.
  
  When $\beta=\alpha+1$ for some $\alpha$, a group in $\LinnelsC_\beta$
  either is an extension of a group in $\LinnelsC_\alpha$ with
  elementary amenable quotient, and thus has the algebraic eigenvalue
  property by Theorem
  \ref{theo:algebraic_eigenvalues_and_amenable_extensions},
  or is a directed union of groups $G_i$ in $\LinnelsC_\alpha$.

  Note that $A\in M(d\times d,KG)$ can be regarded as a matrix $A'$ in
  $M(d\times d,KH)$ where $H$ is a finitely generated subgroup of
  $G$, generated by the finite support of the $A_{i,j}$ in $G$.
  By Proposition 3.1 of \cite{Schick(1998a)}, the spectral
  density functions of $A$ and $A'$ coincide. As subgroups
  of a group with the algebraic eigenvalue property also have
  the property, it follows that a group has the algebraic
  eigenvalue property if and only if it holds for all of its
  finitely generated subgroups.
  If $G\in\LinnelsC_{\alpha+1}$ is the directed union of groups
  $G_i\in \LinnelsC_\alpha$, it follows that every finitely
  generated subgroup of $G$ is in some $G_i$ and so has the
  algebraic eigenvalue property. $G$ therefore has the
  algebraic eigenvalue property.
\end{proof}

In view of these results, we make the following conjecture:
\begin{conjecture}
  Every discrete group $G$ has the algebraic eigenvalue property.
\end{conjecture}
 
We can give further evidence for this conjecture in the case of $G$
belonging to the class $\RAgroups$, using the knowledge of the
spectrum that we have in this case.

\begin{theorem}\label{theo:non_eigenvalue}
  Assume $G\in\RAgroups$ and $A=A^*\in M(d\times d, \overline\rationals G)$.
  Assume that $\lambda\in\reals$ is a transcendental number, but that
  for every $n\in\naturals$ there is a rational number $p_n/q_n$ with
  $q_n\ge 2$ such that
  \begin{equation}\label{eq:good_rational_approxi}
    0< \abs{\lambda-\frac{p_n}{q_n}} \le \frac{1}{q_n^{n}}.
  \end{equation}
  Then $\lambda$ is not an eigenvalue of $A$ acting on $l^2(G)^d$.
\end{theorem}

Observe that it follows from Liouville's theorem \cite[Satz
191]{Hardy-Wright(1958)} that a number $\lambda$ satisfying the second
set of
assumptions of Theorem \ref{theo:non_eigenvalue} is automatically
transcendental. 

\begin{proof}[Proof of Theorem \ref{theo:non_eigenvalue}.]
  The support of the elements $A_{i,j}$ of $A$ over $G$ is finite,
  and so we can find a finite field extension $L\subset \complexs$
  of $\rationals$ and a positive integer $m$ 
  such that $mA\in M(d\times d,o(L)G)$, where $o(L)$ is the ring
  of integers of $L$. If $\lambda$ is a Liouville transcendental,
  then so is $m\lambda$. We can then regard $A$ as being in
  $M(d\times d,o(L)G)$ without loss of generality.

  Since $G\in\RAgroups$, we can apply (\ref{eq:8}) to the operator
  \begin{equation*}
    V_n:=(q_n A- p_n)^*(q_nA-p_n).
  \end{equation*}
  Observe that 
  \begin{equation}
    \label{eq:norm_Vn_estimate}
    \norm{V_n}\le (q_n\norm{A}+\abs{p_n})^2  \le q_n^2 (\norm{A}+\lambda+1)^2,
  \end{equation}
  since by \eqref{eq:good_rational_approxi} $\abs{p_n}\le\abs{\lambda}q_n+q_n$.
  If $\lambda$ is an eigenvalue of $A$, then
  $s_n:=(q_n\lambda-p_n)^2$ will be an eigenvalue of $V_n$ (since
  $A=A^*$). Observe that 
  \begin{equation}\label{eq:s_n_upperbound}
    s_n \le q_n^{2-2n} <1.
  \end{equation}
  If 
  \begin{equation*}
    0< \alpha:= \dim_G(\ker(A-\lambda))
  \end{equation*}
  is the normalized dimension of the eigenspace to $\lambda$, then
  \begin{equation}\label{eq:lower_est_for_eigendim}
    F_{V_n}(s_n)-F_{V_n}(0) = F_{V_n}((q_n\lambda-p_n)^2)-F_{V_n}(0)
    \ge \alpha.
  \end{equation}
  Since $V_n$ has coefficients in $o(L)$, we can use
  (\ref{eq:7}),
  \begin{equation}
    \label{eq:est_expn}
    \begin{split}
      F_{V_n}(s_n)-F_{V_n}(0)
      &\le
      \frac{d\cdot\sum_{k=1}^r \ln(\kappa(\sigma_k(V_n)))}
      {-\ln(s_n/\norm{V_n})}
      \\
      &\le
      \frac{d\cdot\sum_{k=1}^r \ln(\kappa(\sigma_k(V_n)))}{(2n-2)\ln q_n}
    \end{split}
  \end{equation}
  With $A$ self-adjoint,
  \begin{equation}
    \label{eq:kappa_expn}
    \begin{split}
      \kappa(\sigma_k(V_n))
      &=
      \kappa(\sigma_k(q_n^2 A^2 -2p_n q_n A + p_n^2))
      \\
      &=\kappa(q_n^2 \sigma_k(A^2) -2p_n q_n \sigma_k(A) + p_n^2)
      \\
      &=
      S(q_n^2 \sigma_k(A^2) -2p_n q_n \sigma_k(A) + p_n^2)
      \abs{q_n^2 \sigma_k(A^2) -2p_n q_n \sigma_k(A) + p_n^2}_\infty
      \\
      &\leq
      \max\{q_n^2,\abs{2p_n q_n},p_n^2\}\cdot
      \bigl(S(A^2)+S(A)+1\bigr)
      \bigl(\abs{\sigma_k(A^2)}_\infty+\abs{\sigma_k(A)}_\infty+1\bigr)
      \\
      &\leq P(q_n) C_k
    \end{split}
  \end{equation}
  where $P$ is some quadratic polynomial with coefficients depending
  only on $\lambda$, and $C_k$ is a constant
  depending only on $k$ and $A$. Let $C=\max\{C_1,\ldots,C_r\}$.
  Then combining (\ref{eq:lower_est_for_eigendim}),
  (\ref{eq:est_expn}) and (\ref{eq:kappa_expn}),
  \[
  \alpha\leq\frac{dr}{2n-2}
  \cdot\frac{\ln(P(q_n) C)}{\ln q_n}.
  \]
  The right hand side becomes
  arbitrarily small as $n\to\infty$. Consequently, $\alpha=0$,
  i.e.~$\lambda$ is not an eigenvalue of $A$.
\end{proof}

\begin{remark}
  \label{rem:non_eigenvalue}
  It is obvious that Theorem \ref{theo:non_eigenvalue} immediately
  extends to transcendental numbers $\lambda\in\reals$ which have very 
  good approximations by not very ``complex'' algebraic numbers. Here, the 
  complexity of an algebraic number $\xi$ would be measured in terms
  of its denominator (i.e.~how big is $k\in\naturals$ such that $k\xi$ 
  is an algebraic integer), in terms of the degree of its minimal
  polynomial, and in terms of the absolute value of the other zeros of 
  the minimal polynomial.

  Unfortunately, the set of real numbers which admit such
  approximations appears to be of measure zero. At least, this is true by 
  \cite[Satz 198]{Hardy-Wright(1958)} for numbers covered by
  Theorem \ref{theo:non_eigenvalue}, and the proof of Hardy and Wright 
  seems to carry over  without difficulty to the larger set described
  above.
\end{remark}

\section{Zero divisors: from algebraic to arbitrary}
\label{sec:algebraic_to_complex}

In the introduction, we claimed that it suffices to study
$\overline\rationals G$ to decide whether $\complexs G$ satisfies the
zero divisor conjecture. For the readers convenience, we give a proof
of this well known fact here.

\begin{proposition}
  \label{prop:algebraic_to_complex}
  Assume that there are $0\ne a,b\in \complexs G$ with $ab=0$. Then we can
  find $0\ne A,B\in\overline\rationals G$ with $AB=0$.
\end{proposition}
\begin{proof}
  Write $a=\sum_{g\in G}\alpha_g g$, $b=\sum_{g\in G}\beta_g g$. Since 
  only finitely many of the $\alpha_g$ and $\beta_g$ are non-zero,
  they are contained in a finitely generated subfield $L\subset
  \complexs$. We may write
  $L=\overline\rationals(x_1,\dots,x_n)[v]$ with $x_1,\dots,x_n$
  algebraically independent over $\rationals$, and $v$
  algebraic over $\overline\rationals(x_1,\dots,x_n)$ (because of the
  theorem about the primitive element, one such $v$ will do). Upon
  multiplication by suitable non-zero elements of
  $\overline\rationals(x_1,\dots,x_n)$, we may assume that $v$ is
  integral over $\overline\rationals[x_1,\dots,x_n]$,
  i.e.~that there is a polynomial $0\ne p(z)=
  z^n+p_{n-1}z^{n-1}+\dots+p_{0}$ with $p_{j}\in
  \overline\rationals[x_1,\dots,x_n]$, irreducible in
  $\overline\rationals(x_1,\dots,x_n)[z]$, and with $p(v)=0$.
  This can be achieved because 
  the quotient field of the ring of integral elements of $L$ is $L$
  itself.

  Moreover, by multiplying $a$ and $b$ by appropriate non-zero
  elements of $\overline\rationals[x_1,\dots,x_n]$ (the common
  denominator of the $\alpha_g$ or $\beta_g$, respectively), we may
  assume that $0\ne a,b\in
  \overline\rationals[x_1,\dots,x_n][v]G$, with $ab=0$.

  We now proceed to construct a ring homomorphism 
  \begin{equation*}
\phi\colon \overline\rationals[x_1,\dots,x_n][v]\to
\overline\rationals     
\end{equation*}
with $0\ne A:=\phi(a)$ and $0\ne B:=\phi(b)$. Obviously,
$AB=\phi(ab)=0$, and this proves the claim.

To construct $\phi$, observe that substitution of algebraic numbers for 
$x_1,\dots,x_n$ defines a well defined homomorphism
$\phi_0\colon
\overline\rationals[x_1,\dots,x_n]\to\overline\rationals$. Each such
homomorphism can be extended to
$\overline\rationals[x_1,\dots,x_n][v]$, provided
$\phi_0(p(z))$ has a solution in $\overline\rationals$. Since the
highest coefficient of $p$ is $1$, $\phi(p)$ is not a constant
polynomial. Consequently, $\overline\rationals$ being algebraically
closed, the required solution of $\phi_0(p(z))$ and therefore the
extension of the ring homomorphism $\phi_0$ to
$\overline\rationals[x_1,\dots,x_n][v]$ always exists.

Choose now $g, g'\in G$ with $\alpha_{g}\ne 0\ne \beta_{g'}$. Observe 
that $\alpha_g\in\overline\rationals[x_1,\dots,x_n][v]$ is integral
over $\overline\rationals[x_1,\dots,x_n]$ since $v$ is integral. 
The same holds for $\beta_{g'}$. In particular, there exist irreducible
polynomials 
$$
  q(z)=z^m+q_{m-1}z^{m-1}+\dots +q_0,\quad r(z)=z^{m'}+\cdots +r_0
  \in
\overline\rationals[x_1,\dots,x_n][z]
$$
with $q(\alpha_g)=0=r(\beta_{g'})$. Irreducibility implies in
particular that 
\begin{equation*}
q_0=q_0(x_1,\dots,x_n)\ne 0\ne
r_0(x_1,\dots,x_n).
\end{equation*}
By induction on $n$, using the fact that
$\overline\rationals$ is infinite and that any non-zero polynomial
over a field has only finitely many zeros, we show that there are
$d_1,\dots,d_n\in\rationals$ with $q_0(d_1,\dots,d_n)\ne 0\ne
r_0(d_1,\dots,d_n)$.

Consequently, if $\phi_0$ is the corresponding substitution
homomorphism 
\begin{equation*}
\phi_0\colon\overline\rationals[x_1,\dots,x_n]\to
\overline\rationals,
\end{equation*}
$0$ is not one of the zeros of the polynomials
$\phi_0(r)$ or $\phi_0(q)$. Let
\begin{equation*}
\phi\colon\overline\rationals[x_1,\dots,x_n][v]\to\overline\rationals
\end{equation*}
be an extension of $\phi_0$. Then $A_g:=\phi(\alpha_g)$ is a zero of
$\phi(q)=\phi_0(q)$, and $B_{g'}:=\phi(\beta_{g'})$ is a zero of
$\phi(r)$, i.e.~$A_g\ne 0\ne B_{g'}$.

However, with $A=\phi(a)$ and $B=\phi(b)$, $A_g$ is the coefficient of 
$g$ in $A$, and $B_{g'}$ is the coefficient of $g'$ in $B$. It follows 
that $A\ne 0\ne B$, as desired, but $AB=0$.
\end{proof}

\section{Approximation of $L^2$-Betti numbers over the complex\\ group
  ring}\label{sec:complapr}

It is an interesting question, in the situation of Theorem
\ref{newapproxi}, whether the convergence result also holds for matrices
over the complex group ring. We don't know about any
counterexample. Here we will give some positive results.
Let $R$ be a domain.  Recall that $R$ satisfies the \emph{Ore
condition} or equivalently, $R$ is an \emph{Ore domain}, means
that given $r,s \in R$ with $s \ne 0$,
then we can find $r_1,r_2, s_1,s_2 \in R$ with $s_1,s_2 \ne 0$
such that $r_1 s = s_1 r$ and $sr_2 = r s_2$ (some authors call
this the left and right Ore condition).
In this situation we can form the
division ring of fractions whose elements are of the form $s^{-1} r$;
these elements will also be of the form $rs^{-1}$.
\begin{proposition}\label{Clim}
  Assume $G$ is torsion free, $\complexs G$ is an Ore domain and whenever
  $0\ne\alpha\in \complexs G$, $0\ne f\in l^2(G)$ then $\alpha f\ne
  0$. Let $G$ be a subgroup of the direct or inverse limit of a
  directed system of
  groups $G_i$, $i\in I$. If $B\in M(d\times d,\complexs G)$ and $B[i]\in
  M(d\times d, \complexs G_i)$ are constructed as in \ref{newsit}, then
  \begin{equation}\label{blim}
    \dim_G(\ker(B)) = \lim_{i\in I} \dim_{G_i}(\ker(B[i])).
  \end{equation}
\end{proposition}

\begin{corollary}
  The statement of Proposition \ref{Clim} is true if $G$ is torsion
  free elementary amenable or if $G$ is amenable and left-orderable.
\end{corollary}
\begin{proof}
  If $G$ is elementary amenable
  and torsion free then the strong Atiyah conjecture over $\complexs G$ is
  true for $G$ and therefore $\dim_G\ker(\alpha)=0$ if $0\ne\alpha\in
  \complexs G$, hence $\alpha$ has trivial kernel on $l^2(G)$. If $G$ is
  right orderable then $0\ne\alpha\in\complexs G$ has trivial kernel
  on $l^2(G)$ by \cite{Linnell(1992)}. In particular $\complexs G$ has
  no non-trivial zero divisors. Since $G$ is amenable, by an argument
  of Tamari \cite{Tamari(1957)} given in Theorem \ref{tamari} $\complexs G$
  fulfills the Ore
  condition.  Therefore the assumptions of
  Proposition \ref{Clim} are fulfilled.
\end{proof}

For the convenience of the reader we repeat Tamari's result:
\begin{theorem}\label{tamari}
  Suppose $G$ is amenable, $R$ is a division ring (e.g.~a field), and
  $R*G$ is a
  crossed product (e.g.~the group ring $RG$). Assume $R*G$ has no
  non-trivial zero divisors. Then $R*G$ fulfills the Ore condition.
\end{theorem}
\begin{proof}
Suppose we are given $\alpha, \sigma \in R*G$ with
$\sigma \ne 0$.  Without loss of generality, we need to find $\beta,
\tau \in R*G$ with $\tau \ne 0$ such that
$\beta \sigma = \tau \alpha$.  Write $\alpha = \sum_
{g \in G} a_g$ and $\sigma = \sum_{g \in G} s_g$, where $a_g, s_g \in R$ for
all $g \in G$.  Let $Z = \supp \alpha
\cup \supp \sigma$ (where $\supp \alpha$ denotes $\{g \in G \mid a_g \ne 0 \}$,
the support of $\alpha$).  Using the F\o lner condition,
we obtain a finite subset $X$ of $G$ such that
$\sum_{g \in Z} \abs{Xg \setminus X} < \abs{X}$.
Let us write $\beta = \sum_{x \in X} b_x x$ and $\tau = \sum_{x \in
X} t_x x$ where $b_x,t_x \in R$ are to be determined.  We want to solve
the equation $\beta \sigma = \tau \alpha$, which when written out
in full becomes
\[
\sum_{x \in X, \ g \in G} b_x(x s_g x^{-1}) xg =
\sum_{x \in X, \ g \in G} t_x(x a_g x^{-1}) xg.
\]
By equating coefficients, this yields at most $2 \abs{X} -1$ homogeneous
equations in the $2\abs{X}$ unknowns $b_x,
t_x$.  We deduce that there exist $\beta, \tau$, not both zero, such that
$\beta \sigma = \tau \alpha$.  Since $R*G$ is a domain and $\sigma
\ne 0$, we see that $\tau \ne 0$ is not a zero divisor and the result
is proved.
\end{proof}

Proposition \ref{prop:resapproxi} is a direct consequence of the corollary
because there $G$ is a subgroup of the inverse limit of the quotients $G/G_i$.

\begin{proof}[Proof of Proposition \ref{Clim}]
  By \cite[7.3]{Schick(1998a)} \eqref{blim} holds if $d=1$. We use the
  Ore condition to reduce the case of matrices to the case
  $d=1$. Since $\complexs G\subset l^2(G)$ the assumptions imply that
  $\complexs G$ does not have zero divisors. Therefore the Ore
  localization $DG$ of $\complexs G$ is a skew field. Moreover, since
  $0\ne\alpha\in\complexs G$ has no kernel on $l^2(G)$ it becomes
  invertible in the ring $\universalU G$ of operators affiliated to
  the group von Neumann algebra of $G$ \cite[\S 8]{Linnell(1998)} and
\cite[Definition 7, p.~741]{Schick(1999)}.
    Therefore $DG$ embeds into $\universalU G$. (By
  \cite[4.4]{Schick(1999)} $G$ fulfills the strong Atiyah conjecture
  over $\complexs G$.)
  
  Fix now $B\in M(d\times d,\complexs G)$. By linear algebra, we find
  invertible matrices $X,Y\in M(d\times d,DG)$ such that
  $XBY=\diag(v_1,\dots,v_d)$ with $v_i\in DG$. The Ore condition
  implies that we can find $a\in \complexs G-\{0\}$
  such that $v_i=\alpha_i a^{-1}$ with $\alpha_i\in\complexs G$,
for $i = 1, \dots, d$. Therefore
  $XBYa=\diag(\alpha_1,\dots,\alpha_d)$
Applying the same principle
  to the entries of $X$ and $Ya$, we can find $x,y\in\complexs
  G-\{0\}$ such that $xX = U$ and $Yay=V$ with $U,V\in
  M(d\times d,\complexs G)$. Altogether we arrive at
  \begin{equation}\label{diagonalize}
      UBV = \diag(a_1,\dots,a_d)
    \end{equation}
    where all objects are defined over $\complexs G$
    ($a_k=x\alpha_ky$). Moreover, 
    $U$ and $V$ are invertible over $DG$ and therefore also over
    $\universalU G$ and hence have trivial kernel as operators on
    $l^2(G)^d$. \eqref{diagonalize} translates to 
    \begin{equation*}
      U[i] B[i] V[i] = \diag((a_1)_i,\dots,(a_d)_i) 
    \end{equation*}
(if $G$ is the direct limit of the groups $G_i$, then the images of
the left and right hand sides of the above in $M(d\times
d, \mathbb {C}G)$ are equal, consequently for all sufficiently large
$i$ the above equality will be true).

    By \cite[7.2]{Schick(1998a)} $\dim_{G_i}\ker(U[i])\xrightarrow{i\in
      I} 0$ and $\dim_{G_i}\ker(V[i])\xrightarrow{i\in I} 0$. The one
    dimensional case immediately implies
    \begin{equation*}
\dim_{G_i}\ker(\diag({a_1}_i,\dots,{a_d}_i))
    \to \dim_G\ker(\diag(a_1,\dots,a_d)).
  \end{equation*} 
    We have the exact sequences
    \begin{gather}
      0\to \ker(V[i])\into \ker(U[i]B[i]V[i]) \xrightarrow{V[i]}
      \ker(U[i]B[i])\\
      0\to \ker(B[i])\into \ker(U[i]B[i])\xrightarrow{B[i]} 
\ker(U[i]).
    \end{gather}
Because of
    additivity of the $L^2$-dimension \cite[Lemma
1.4(4)]{Lueck(2001)}
    \begin{multline*}
     \dim_{G_i}\ker(U[i]B[i]V[i])-\dim_{G_i}\ker(U[i])-\dim_{G_i}\ker(V[i])\\
     \le \dim_{G_i}\ker(B[i])\le \dim_{G_i}\ker(U[i]B[i]).
   \end{multline*}
   Since all these operators are endomorphism of the same finite
   Hilbert $\NeumannN G_i$-module $l^2(G_i)^d$, 
   \begin{multline*}
\dim_{G_i}
\ker(U[i]B[i])=\dim_{G_i}\ker(B[i]^*U[i]^*)\\
  \le
   \dim_{G_i}\ker(V[i]^*B[i]^*U[i]^*)=\dim_{G_i}\ker(U[i]B[i]V[i]), 
 \end{multline*}
   and similarly $\dim_G\ker(UBV)=\dim_G\ker(B)$. Everything together
   implies \eqref{blim}.
\end{proof}

\bibliographystyle{plain}

\end{document}